%% file: 2005-43.tex
\documentclass{gtart_h}  

\input gtoutput

\lognumber{497}
\received{30 January 2004}
\volumenumber{9}\papernumber{43}\volumeyear{2005}
\pagenumbers{1915}{1951}   
\revised{6 October 2005}
\published{17 October 2005}
\accepted{30 August 2005}
\proposed{Benson Farb}
\seconded{Martin Bridson, Walter Neumann}

\usepackage{amsmath,amssymb,graphicx,psfrag}

\def\psfraga <#1,#2> #3#4{%
\psfrag {#3}{\smash{\rlap{\kern #1 \raise #2\hbox{#4}}}}}

\def\figref#1{\hyperlink{#1anchor}{Figure~\ref*{#1}}}
\def\anchor#1{\noindent\hypertarget{#1anchor}{\smash{$\phantom{99}$}}\newline}

\let\tilde\widetilde

\theoremstyle{definition}
\newtheorem{dfn}{Definition}[section]
\newtheorem{rem}[dfn]{Remark}

\newtheorem{defn}[dfn]{Definition}
\newtheorem{definition}[dfn]{Definition}
\newtheorem{problem}[dfn]{Problem}

\theoremstyle{plain}
\newtheorem{thm}[dfn]{Theorem}
\newtheorem{lem}[dfn]{Lemma}
\newtheorem{prop}[dfn]{Proposition}
\newtheorem{obs}[dfn]{Observation}
\newtheorem{cor}[dfn]{Corollary}
\newtheorem{conj}[dfn]{Conjecture}

\def\no{\noindent}

\newcommand{\restr}{\mbox{\Large \(|\)\normalsize}}
\def\embed{\hookrightarrow}

\def\bul{\bullet}
\def\ol{\overline}
\def\0{\emptyset}

\def\>{\rangle}
\def\<{\langle}

\def\C{\mathbb C}
\def\t{\tilde}
\def\R{\mathbb R}

\def\p{{\cal P}}

\def\b{{\mathcal B}}
\def\e{{\mathcal E}}

\def\Q{\mathbb Q}

\def\A{{\mathbb A}}

\def\H{\mathbb H}

\def\al{\alpha}

\def\ga{\gamma}
\def\Ga{\Gamma}

\def\Si{\Sigma}
\def\si{\sigma}

\def\La{\Lambda}
\def\8{\infty}

\def\D{\partial}

\def\geo{\D_\8 }

\def\1{\hbox{\bf 1}}

\def\ov{\overrightarrow}
\def\bul{\bullet}

\def\p{{\mathcal P}}

\def\Length{{\rm Length}}
\def\Span{{\rm Span}}
\def\Bis{{\rm Bis}}

\newcommand{\isom}{\operatorname{Isom}}

\newcommand{\im}{\operatorname{Im}}
\newcommand{\acts}{\curvearrowright}
\renewcommand{\k}{{\bf k}}

\begin{document}
\title{Representations of polygons of finite groups}
\author{Michael Kapovich}
\address{Department of Mathematics, University of California\\
1 Shields Ave, Davis, CA 95616, USA}
\email{kapovich@math.ucdavis.edu}

\begin{abstract}
We construct discrete and faithful representations into the isometry
group of a hyperbolic space of the fundamental groups of acute
negatively curved even-sided polygons of finite groups.
\end{abstract}

\primaryclass{20F69} 
\secondaryclass{20G05}
\keywords{Hyperbolic groups, polygons of groups}
\maketitlepage
%\makeagttitle

\section{Introduction}
There are only few known obstructions for existence of an isometric properly discontinuous 
action of a Gromov-hyperbolic group $G$ on the {\em real-hyperbolic space} $\H^p$ for some $p$:

(1)\qua If $G$ is a group satisfying Kazhdan property (T) then each isometric action $G\acts \H^p$ 
fixes a point in  $\H^p$; hence no infinite hyperbolic group satisfying property (T) admits 
 an isometric properly discontinuous action $G\acts \H^p$, for {\em any $p$}. 
 
(2)\qua Suppose that $G$ is the fundamental group of a compact K\"ahler manifold and 
$G\acts \H^p$ is  an isometric properly discontinuous action. Then, according to 
a theorem of Carlson and Toledo \cite{Carlson-Toledo}, this action factors 
through an epimorphism $G\to Q$, where $Q$ is commensurable to a surface group. Hence, 
unless $G$ itself is commensurable to a surface group, it does not admit an 
 isometric properly discontinuous action $G\acts \H^p$. Examples of Gromov-hyperbolic groups 
 which are K\"ahler (and are not commensurable to surface groups) 
 are given by the uniform lattices in $PU(m, 1)$, $m\ge 2$, as well as the fundamental 
 groups of compact negatively curved K\"ahler manifolds (see \cite{Mostow-Siu}).  

On the positive side, by a theorem of Bonk and Schramm \cite{bonkschramm}, 
each Gromov-hyperbolic group admits a quasi-isometric embedding to a real-hyperbolic space.

The goal of this paper is to find a better ``demarcation line'' between hyperbolic groups 
satisfying property (T) and groups acting discretely on real-hyperbolic spaces. 
In this paper we will show that a large class of 2--dimensional Gromov-hyperbolic groups 
admits isometric properly discontinuous {\em convex-cocompact} actions on real-hyperbolic spaces. 
We consider a 2--dimensional negatively curved acute polygon $\p$ of finite groups (see 
section \ref{polygons} for more details). 
Let $G:= \pi_1(\p)$ be the fundamental group of this polygon, we refer the reader 
to \cite[Chapter II, section 12]{Bridson-Haefliger} for the precise definitions. 
%Note that negative curvature implies that the polygon of groups $\p$ is {\em developable}  
%(see \cite{Bridson-Haefliger});  
%let $G\acts X$ denote the corresponding cocompact isometric action 
%on the universal cover $X$ of $\p$. 

Our main result is:

\begin{thm}\label{main}
Suppose that $n=2k$ is even. Then the group $G$ admits a discrete, 
faithful, convex-cocompact action $\rho$ on a constant curvature hyperbolic space $\H^p$, 
where $p<\infty$ depends on the polygon $\p$.  
\end{thm}   

Our technique in general does not work in the case when $n$ is odd: 
We were unable to construct 
a representation. However in section \ref{odd} we will construct 
$\rho$ and prove that it is discrete, 
faithful, convex-cocompact for a special class of odd-sided 
$n$--gons of groups, provided that $n\ge 5$. 

\begin{conj}
The assertion of Theorem \ref{main} remains valid for all odd $n\ge 5$.  
\end{conj}

\medskip 
In contrast, if $\p$ is a %negatively curved 
{\em triangle} of finite groups where the vertex links are connected graphs 
with the $1^{\rm st}$ positive eigenvalue of the Laplacian $ >1/2$, 
then the group $G=\pi_1(\p)$ satisfies property (T), see \cite{Ballmann-S(1997)}.  
Hence (provided that $G$ is infinite) the group $G$ cannot act properly discontinuously on 
$\H^p$ for any $p$. Thus, it appears, that (at least for the polygons of finite groups) 
the ``demarcation line'' which we are trying to find, is hidden somewhere among quadrilaterals   
and triangles of groups. We will address this issue in another paper.  

\medskip 
Recall that, by a theorem of Dani Wise \cite{Wise}, $G$ is residually finite 
(actually, Wise proves that $G$ has separable quasi-convex subgroups, 
which is used in the proof of our main theorem). 

\begin{cor}
The group $G$ is linear. 
\end{cor}

\begin{rem}
A very different proof of linearity of $G$ was given by Wise and Haglund, who used an embedding 
of $G$ to a right-angled Coxeter group, \cite{HW}.  
\end{rem}

The following problem is open even for right-angled Coxeter groups of virtual 
cohomological dimension 2.

\begin{problem}
Suppose that $G$ is a Gromov-hyperbolic Coxeter group. Is $G$ isomorphic to a 
discrete subgroup of $\isom(\H^n)$ for some $n$? Note that if one insists that 
the Coxeter generators act on $\H^n$ as reflections, then there are examples 
of  Gromov-hyperbolic Coxeter  groups which do not admit such actions on $\H^n$, see 
\cite{Felikson-Tumarkin}. 
\end{problem}

\medskip 
In  section \ref{appendix}, we give an example of a nonlinear Gromov-hyperbolic group.

%On the other hand, our construction in the even case works--to some extent 
%-- for $n=4$: We get both discrete (but non-faithful) 
%representations into $\isom(\R^p)$ and (in general non-discrete) 
%representations into $\isom(\H^p)$, which do not have a fixed 
%point (even at infinity). 

The proof of the main theorem splits in two parts: (1) Construction of $\rho$, 
(2) proof of discreteness. To prove discreteness of $\rho$ we show that there 
exists a $\rho$--equivariant quasi-isometric 
embedding $\mu\co  X\to \H^p$, where $X$ is  the universal cover of the polygon $\p$. 
This proves that the action $\rho\co  G\acts \H^p$ is properly 
discontinuous and convex-cocompact. (A priori this action can have finite kernel. In 
section \ref{extension} we explain how to deal with this issue.)  
The proof that $\mu$ is a quasi-isometric embedding 
is based on the following theorem of independent interest:

\begin{thm}
%\label{qi}
Suppose that $X$ is a 2--dimensional regular cell complex, 
which is equipped with a ${\rm CAT}(-1)$ path-metric so that 
each face of $X$ is isometric to a right-angled regular $n$--gon in $\H^2$. 
Let $\mu\co  X\to \H^p$ is a continuous map which is a (totally-geodesic) isometric embedding 
on each face of $X$. Assume also that for each pair of faces $F', F''\subset X$ which 
intersect non-trivially a common face $F\subset X$, we have:
$$
\Span(\mu(F')) \perp \Span(\mu(F'')). 
$$
Then $\mu$ is a quasi-isometric embedding. 
\end{thm}

Our construction of representations $\rho$ was inspired by the paper of Marc Bourdon, 
\cite{Bourdon97}, 
where he proves a theorem which is a special case of Theorem \ref{main}: 
In his paper Bourdon considers  
$n$--gons of finite groups where the edge groups are cyclic and the 
vertex groups are direct products of the adjacent edge groups, 
under the extra assumption that the orders of the edge groups are much smaller that $n$.  
%To prove discreteness of his representation Bourdon uses Dirichlet fundamental domains. 
%Note that the first part of the proof does not require $\mu(F)$ to be regular, nor right-angled, although we need 
%all angles to be the same. This leads to nontrivial deformations of representations constructed in the paper. 
%However, if one side of $\mu(F)$ is allowed to shrink, the representation eventually will fail to be
%discrete (or faithful).  

\medskip 
{\bf Acknowledgments}\qua During the work on this paper I was visiting
the Max Plank Institute (Bonn), I was also supported by the NSF grants
DMS-02-03045 and DMS-04-05180.  I am grateful to Tadeusz Januszkiewicz
for inspirational and helpful conversations during the work on this
paper. I first tried to prove discreteness of certain representations
by verifying that an equivariant map is a quasi-isometric embedding,
in a joint project with Bernhard Leeb in 1998. Although our attempt
back then was unsuccessful, I am grateful to Bernhard Leeb for that
effort.  I am also grateful to the referees of this paper for numerous
suggestions and to Mark Sapir for discussions of Theorem \ref{non}.

%\tableofcontents 

\section{Preliminaries}

{\bf Notation}\qua If $\Si$ is a (finite) set, we define a Euclidean vector 
space $Vect(S)$ to be the vector space $L_2(S)$, where $S$ forms an orthonormal 
basis (we identify 
each 1--point subset of $S$ with its characteristic function in $L_2(S)$). 
Suppose that $S\subset \H^q$. Then $\Span(E)$ will denote the smallest 
totally-geodesic subspace in $\H^q$ which contains $S$. If $E\subset \H^q$ 
is a geodesic segment, then $\Bis(E)$ will denote the perpendicular 
bisector of $E$. 

Suppose that $E_1, E_2\subset E$ are subspaces of a Euclidean vector space $E$, whose 
intersection is $E_3$. 
We say that $E_1, E_2$ {\em intersect orthogonally} if $E_1/E_3, E_2/E_3$ 
are contained in the orthogonal complements of each other in the Euclidean space $E/E_3$. 

Suppose that $H', H''$ are totally-geodesic subspaces in $\H^p$. 
We say that $H', H''$ {\em intersect orthogonally} if 
$$
H'\setminus H''\ne \emptyset, H''\setminus H'\ne \emptyset
$$ 
and for some (equivalently, for every) point $x\in H'\cap H''$ we have:   
$$T_x(H'), T_x(H'')\subset T_x(\H^p) \hbox{~~intersect orthogonally}.$$ 

Totally geodesic subspaces $H', H''\subset \H^p$ are said to be 
{\em orthogonal to each other} if either:

(a)\qua $H', H''$ intersect orthogonally, or 

(b)\qua $H', H''$ are within positive distance from each other and for the unique shortest 
geodesic segment $\si:=\ol{x'x''}\subset \H^p$ connecting $H'$ to $H''$, the totally-geodesic 
subspaces $H'', \ga_{\si}(H')$ intersect orthogonally. Here $\ga_{\si}$ is 
the hyperbolic translation along $\si$ which sends $x'$ to $x''$. 

\medskip 
We will use the notation $H'\perp H''$ for subspaces $H', H''$ orthogonal to each other. 
Clearly, $H'\perp H'' \iff H''\perp H'$.

\subsection[Discrete subgroups of isom(H^n)]{Discrete subgroups of $\isom(\H^n)$}

Recall that a map $f\co  X\to Y$ between two metric spaces is called 
an $(L,A)$ {\em quasi-isometric embedding} if for all $x, x'\in X$ we have:
$$
L^{-1}d(x, x') -A \le d(f(x), f(x'))\le L d(x, x') +A,
$$
where $L>0$. 
An $(L,A)$ {\em quasi-isometry} is an $(L,A)$ quasi-isometric embedding 
$X \stackrel{f}{\to} Y$ such that each point of $Y$ is within distance $\le A$ 
from a point in $\im(f)$.  

A map $f$ is called a  {\em quasi-isometry} (resp.\ a {\em quasi-isometric embedding})  
if it is an  $(L,A)$ quasi-isometry (resp.\ quasi-isometric embedding) 
for some $L$ and $A$. 

An $(L,A)$ quasi-geodesic segment in a metric space $X$ is an 
$(L,A)$ quasi-isometric embedding $f\co  [0, T]\to X$, where $[0, T]$ is 
an interval in $\R$. By abusing notation we will sometimes refer 
to the image $\im(f)$ of an $(L,A)$ quasi-geodesic segment $f$ as an 
$(L,A)$ quasi-geodesic segment.  
Recall that by the {\em Morse lemma} (see for instance \cite[Lemma 3.43]{Kapovich2000}), 
quasi-geodesics in $\H^n$ are {\em stable}: 

There is a function $D=D(L,A)$ such that for each 
$(L,A)$ quasi-geodesic segment $f\co  [0, T]\to \H^n$, 
the Hausdorff distance between $\im(f)$ and the geodesic segment 
$$
\ol{f(0) f(T)}\subset \H^n
$$ 
connecting the end-points of $\im(f)$ is at most $D$.  

\medskip 
We will use the notation $\H^n$ for the real-hyperbolic $n$--space; 
its curvature is normalized to be equal to $-1$. 
The space $\H^n$ has a geometric compactification $\bar{\H}^n= \H^n \cup S^{n-1}$. 
For a subset $S\subset \H^n$ we 
let $\bar{S}$ denote its closure in $\bar{\H}^n$. 

Discrete subgroups $G\subset \isom(\H^n)$ are called {\em Kleinian groups}. 

The {\em convex hull} $C(G)$ of 
a Kleinian group  $G\subset \isom(\H^n)$ is the smallest  nonempty closed 
convex $G$--invariant subset $C\subset \H^n$.  
The convex hull exists for each $G$ whose limit set has cardinality $\ne 1$. 
The convex hull is unique unless $G$ is finite.

\begin{defn}
A Kleinian group $G\subset \isom(\H^n)$ is called {\em convex-cocompact} 
if $C(G)$ exists and the quotient $C(G)/G$ is compact. 
\end{defn}
 
\begin{lem}
\label{qconvex}
Suppose we have a representation $\rho\co G\to \isom(\H^n)$ of a
finitely-generated group $G$. Then the action $\rho\co G\acts \H^n$ is
properly discontinuous and convex-cocompact iff there exists a
$G$--equivariant quasi-isometric embedding $f\co \Ga_G\to \H^n$, where
$\Ga_G$ is a Cayley graph of $G$.
\end{lem} 
\proof First, suppose that $\rho\co  G\acts \H^n$ properly discontinuous 
and convex-cocompact. Then, 
because $C(G)$ is a geodesic metric space, there exists  
a $G$--equivariant quasi-isometry $f\co  \Ga_G\to C(G)$. Composing this map with the 
isometric embedding $\iota\co  C(G)\to \H^n$, we conclude that 
$f\co  \Ga_G\to \H^n$ is a quasi-isometric embedding. 

Conversely, suppose that $f\co  \Ga_G\to \H^n$ is an equivariant 
quasi-isometric embedding. In particular, $f$ is a proper map. 
Hence, if for $1\in \Ga_G$ we set $o:=f(1)$, then for each compact subset 
$K\subset \H^n$ there are only finitely many elements $g\in G$ such 
that $g(o)\in K$. Therefore the action $G\acts \H^n$ is properly 
discontinuous. In particular, it has finite kernel. 

Observe that {\em stability} of quasi-geodesics in $\H^n$ implies that $\im(f)$ 
is {\em quasi-convex}, ie, there exists a constant $c<\infty$ such that 
for any two points $x, y\in \im(f)$ the geodesic segment $\ol{xy}$ is contained 
in a $c$--neighborhood $N_c(\im(f))$ of $\im(f)$. 
On the other hand, by \cite[Proposition 2.5.4]{Bowditch(1995)}, 
there exists $R=R(c)$ such that the convex hull of each $c$--quasi-convex subset
$S\subset \H^n$ is contained in the $R$--neighborhood $N_R(S)$. Thus, the 
convex subset $C(\rho(G))\subset \H^n$ is contained 
in $N_{R(c)}(\im(f))$. Since $G$ acts cocompactly on $\im(f)$ it follows 
that $G$ acts cocompactly on $C(G)$. Therefore $\rho(G)$ is convex-cocompact. 
\qed 

\begin{rem}
(1)\qua Clearly, instead of $\Ga_G$ in the above lemma 
one can use any geodesic metric space on which $G$ 
acts isometrically, properly discontinuously and cocompactly.

(2)\qua The above lemma shows that existence of a $G$--equivariant quasi-isometric embedding   
$f\co  \Ga_G\to \H^n$ implies that $\rho$ has finite kernel. 
However it does not exclude the possibility that this kernel is nontrivial.  
\end{rem}

\begin{defn}
Let $G$ be a group with a Cayley graph $\Ga_G$. A subgroup $H\subset G$ is called 
{\em quasi-convex} if the orbit $H\cdot 1\subset \Ga_G$ is quasi-convex, ie, 
there exists a number $D$ so that each geodesic segment 
$\si\subset \Ga_G$ with vertices in $H\cdot 1$ is contained in $N_D(H\cdot 1)$.  
\end{defn}

If $G$ is Gromov-hyperbolic then quasi-convexity of $H$ 
is independent of the choice of Cayley graph $\Ga_G$. 

\subsection{Geometry of polygons of groups}
\label{polygons}

Consider an $n$--gon $P$ ($n\ge 5$) with vertices $x_i$ and edges $e_j$, $1\le i, j\le n$. 
Throughout we will be working mod $n$, ie, $qn+i$ will be identified with $i$ 
for $i\in \{1,\ldots, n\}$. We will be assuming that each edge $e_i$ has 
the vertices $x_i, x_{i+1}$. We will regard $P$ as a (2--dimensional) cell complex 
and its poset $Pos(P)$ as a (small) category. A {\em polygon of groups} $\p$ based on $P$ 
is a covariant functor from $Pos(P)$ to the category of groups and monomorphisms.  

In other words, a {\em polygon of groups} $\p$ based on $P$ is a collection 
of groups $G_{x_i}, G_{e_i}, G_F$ assigned to the vertices, edges and the 2--face $F$ 
of $P$, together with monomorphisms 
$$
G_F\to G_e\to G_x$$
for each edge $e$ containing the vertex $x$, so that the following diagrams are commutative: 
$$
\begin{array}{ccc}
G_x& \leftarrow & G_e\\
~ & \nwarrow & \uparrow\\
~ & ~ & G_F
\end{array}
$$   
The direct limit of the above diagrams of monomorphisms is the {\em fundamental group}  
$G=\pi(\p)$ of the polygon $\p$. If the vertex, edge and face groups of $\p$ 
embed naturally into $G$, the polygon $\p$ is called {\em developable}. 
Not every polygon of groups is developable, however under a certain {\em nonpositive curvature} 
assumption on $\p$, the polygon $\p$ is developable, see \cite{Bridson-Haefliger}. 

\medskip 
{\bf Curvature and angles}\qua For each vertex $x_i\in P$ define a graph $Lk_{x_i}$ as follows. 
The vertices of $Lk_{x_i}$ are the right cosets $gG_{e_{i}}, gG_{e_{i-1}}, g\in G_{x_i}$. 
The vertices $v, w$ are connected by a (single) edge iff there exists $g\in G_{x_i}$ 
such that $g(\{v, w\})= \{G_{e_{i}}, G_{e_{i-1}}\}$. Thus the group $G_{x_i}$ acts 
on $Lk_{x_i}$ with the quotient being the edge connecting 
$G_{e_{i}}, G_{e_{i-1}}$. We metrize the graph $Lk_{x_i}$ by assigning the same 
length $\al_i$ to each edge, so that the group $G_{x_i}$ acts isometrically. 
Then the {\em angle} between the subgroups $G_{e_{i}}, G_{e_{i-1}}$ is the least 
number $\al_i$ such that the metric graph $Lk_{x_i}$ is a CAT(1) space, ie, 
the length of the shortest embedded cycle in $Lk_{x_i}$ is at least $2\pi$. 
Equivalently, the angle between $G_{e_{i}}, G_{e_{i-1}}$ equals
$$
2\pi/girth(Lk_{x_i}). 
$$
We will say that the polygon $\p$ is {\em acute} (or has {\em acute angles}) 
if the angle $\al_i$ between each pair of edge groups $G_{e_{i}}, G_{e_{i-1}}$ is at most $\pi/2$.  

We refer the reader to 
\cite[Chapter II, section 12]{Bridson-Haefliger}  
for the precise definitions of the nonpositive/negative curvature of $\p$; various examples 
of negatively curved polygons of groups can be found in \cite{Ballmann-S(1997)}, 
\cite[Chapter II, section 12]{Bridson-Haefliger} and \cite{Swiatkowski(2001)}. 
Instead, we state the following equivalent definition of negative curvature:  
%corollary of negative curvature for acute polygons of groups: 

There exists a 2--dimensional simply-connected regular cell complex $X$ 
(the universal cover of $\p$) together with a path-metric on $X$ whose restriction 
to each face of $X$ has constant curvature $-1$, so that: 

\medskip 
(1)\qua Each face of $X$ is isometric to an $n$--gon in $\H^2$ with angles $\al_1,\ldots,\al_n$. 

(2)\qua Each cell in $X$ is convex. 

(3)\qua There exists an isometric cellular action $G\acts X$ which is transitive on 2--cells.

(4)\qua The stabilizer of each 2--face $F\subset X$ is isomorphic to $G_F$, it fixes $F$ pointwise.

(5)\qua The stabilizer of each edge $e$ of $F$ is isomorphic to $G_e$ and it fixes $e$ pointwise.

(6)\qua The stabilizer of each vertex $x$ of $F$ is isomorphic to $G_x$.

(7)\qua The inclusion maps $G_F\embed G_e\embed G_x$ coincide with the monomorphisms 
$G_F\to G_e\to G_x$ in the definition of $\p$.  

\medskip
Note that the link in $X$ of each vertex $x_i\in F$ is isometric to $Lk_{x_i}$ 
(where each edge has the length $\al_i$). Thus the above complex $X$ is a CAT(-1) metric space. 

Throughout the paper we will be using only the following corollary of negative curvature 
for {\em acute} polygons of groups:

\begin{cor}
If $\p$ is negatively curved then there is a ${\rm CAT}(-1)$ complex $X$ where each face is 
isometric to a regular right-angled polygon in $\H^2$, so that the properties (2)--(7) are 
satisfied. 
\end{cor}

\smallskip 
In this paper we will consider only the case when the vertex groups are finite, 
thus the action $G\acts X$ is properly discontinuous and cocompact, which 
implies that $X$ is equivariantly quasi-isometric to a Cayley graph of $G$. 

We now return to the original polygon $\p$ assuming that it has even number of sides. 
Let $o$ denote the center of the face $F$ and let $m_j$ be the midpoint of the edge $e_j\subset F$.  
We consider two subgraphs $\Ga_{even}, \Ga_{odd}\subset F$ 
which are obtained by conning off from $o$ the sets 
$$
m_{even}:=\{m_{2j}, j=1,\ldots, n/2\}, m_{odd}:=\{m_{2j-1}, j=1,\ldots, n/2\}
$$
respectively. 
Let $G_{even}, G_{odd}$ denote the subgroups of $G$ generated by the elements of 
$$
G_{e_{2j}}, j=1,\ldots,n/2
$$
and
$$
G_{e_{2j-1}}, j=1,\ldots,n/2
$$
respectively. Define subgraphs $T_{even}$ and  $T_{odd}$ to be the orbits 
$$
G_{even}\cdot \Ga_{even}\quad \hbox{and}\quad G_{odd}\cdot \Ga_{odd}.
$$
 We define a new path-metric $\tau$ on the complex $X$ by declaring the closure 
of each component of $X\setminus (T_{even}\cup T_{odd})$ to be a unit Euclidean square. 
Clearly, the group $G$ acts on $(X, \tau)$ isometrically and  $(X, \tau)$ is a CAT(0) metric space.

The groups $G_{even}, G_{odd}$ act on the graphs  $T_{even}, T_{odd}$ with the 
fundamental domains $\Ga_{even}, \Ga_{odd}$ respectively. 
It therefore follows that if $g\in G$ and $g(T_{even})\cap T_{even}\ne \emptyset$ 
 (resp.\  $g(T_{odd})\cap T_{odd}\ne \emptyset$) then $g\in G_{even}$ (resp.\ $g\in G_{odd}$). 

\begin{lem}
The subgraphs $T_{even}, T_{odd}\subset X$ are convex subsets in $X$ isometric 
(with respect to the path-metric induced from $(X,\tau)$) to a tree. 
\end{lem}
\proof (1)\qua First, let us prove that $T_{even}, T_{odd}\subset X$ are convex. Since $X$ is a 
{\rm CAT(0)} space, and $T_{even}, T_{odd}$ are connected, it suffices to test convexity at each 
vertex of $T_{even}, T_{odd}$. However, by the definition of the metric $\tau$, the angle 
between different edges of $T_{even}$ (resp.\ $T_{odd}$) at each vertex of $T_{even}$ 
(resp.\ $T_{odd}$) is $\ge \pi$. Therefore convexity follows. 

(2)\qua Since, $T_{even}, T_{odd}\subset (X, \tau)$ are convex, it follows that they are contractible. Therefore 
these graphs are isometric to metric trees. It is clear that $T_{even}$ and $T_{odd}$ 
are isometric to each other. \qed  
 
\begin{cor}
Each subgroup $G_{even}, G_{odd}$ is a quasi-convex subgroup of $G$. 
\end{cor}
\proof Cayley graphs $\Ga_G$, $\Ga_{G_{even}}$, $\Ga_{G_{odd}}$ are quasi-isometric 
to $(X, \tau)$, $T_{even}$, $T_{odd}$ respectively. Recall that $X$ is Gromov-hyperbolic. 
Therefore, by combining stability of quasi-geodesics in Gromov-hyperbolic geodesic metric spaces 
and convexity of $T_{even}, T_{odd}\subset X$, we conclude that 
 $G_{even}, G_{odd}$ are quasi-convex subgroups of $G$. \qed  

\medskip 
We define functions $odd(i)$ and $even(i)$ by
$$
odd(i)=
\left\{\begin{array}{c}
i, \hbox{~if $i$ is odd}\\
i-1,  \hbox{~if $i$ is even}
\end{array}\right. \hbox{~and~}\quad  
even(i)=
\left\{\begin{array}{c}
i, \hbox{~if $i$ is even}\\
i-1,  \hbox{~if $i$ is odd}
\end{array}\right.
$$
To motivate the following definition, observe that the group $G$ is generated by the elements 
$g_l$ of the vertex subgroups $G_{x_l}$. It will be very important for the later analysis 
to find out which products of pairs of generators $f=h_j^{-1} g_i, h_j\in G_{x_j}, g_i\in G_{x_i}$,  
preserve the trees $T_{odd}, T_{even}$. Note that the answer is clear for some of these  
products:  

(a)\qua If $g_i$ preserves $e_{odd(i)}$ (resp.\ $e_{even(i)}$) and $h^{-1}_j$ preserves $e_{odd(j)}$ 
(resp.\ $e_{even(j)}$), then $g_i, h_j^{-1}\in G_{odd}$ (resp.\ $G_{even}$) and hence 
$f$ also preserves  $T_{odd}$ (resp.\ $T_{even}$). 

(b)\qua If $i=j$ and the product $f=h_i^{-1} g_i$ preserves the edge $e_{odd(i)}$ (resp.\ $e_{even(i)}$), 
then $f$ also preserves  $T_{odd}$ (resp.\ $T_{even}$).

\medskip
\no Accordingly, define finite subsets $\Phi_{even}', \Phi_{even}, \Phi_{odd}', 
\Phi_{odd}\subset G$ as follows:

(1)\qua $\Phi_{even}'$ consists of products $h_j^{-1} g_i$, $g_i\in G_{x_i}, h_j\in G_{x_j}$ with either 
(a) $i\ne j$ and $g_i\in G_{e_{even(i)}}$ and $h_j\in G_{e_{even(j)}}$, or (b) $i=j$ and 
$h_i^{-1} g_i\in G_{e_{even(i)}}$.   

(2)\qua $\Phi_{even}:= \{ h_j^{-1} g_i: g_i\in G_{x_i}, h_j\in G_{x_j}\} \setminus \Phi_{even}'$. 

(3)\qua  $\Phi_{odd}'$ consists of products $h_j^{-1} g_i$, $g_i\in G_{x_i}, h_j\in G_{x_j}$ with either 
(a) $i\ne j$ and $g_i\in G_{e_{odd(i)}}$ and $h_j\in G_{e_{odd(j)}}$, or (b) $i=j$ and 
$h_i^{-1} g_i\in G_{e_{odd(i)}}$. 

(4)\qua $\Phi_{odd}:= \{ h_j^{-1} g_i: g_i\in G_{x_i}, h_j\in G_{x_j}\} \setminus \Phi_{odd}'$.

\no Observe that $\Phi_{even}', \Phi_{odd}'$ are contained in the subgroups $G_{even}, G_{odd}$ 
respectively. 

\begin{prop}
$\Phi_{odd}\cap G_{odd}=\emptyset$, $\Phi_{even}\cap G_{even}=\emptyset$.  In other words, 
among the products of the generators, only the ``obvious'' ones preserve the trees 
$T_{odd}$ and $T_{even}$. 
\end{prop}
\proof We prove that $\Phi_{odd}\cap G_{odd}=\emptyset$, the second assertion is 
proved by relabeling. We have to show that 
$$
 h_j^{-1} g_i T_{odd}= T_{odd} \Rightarrow h_j^{-1} g_i\in \Phi_{odd}'. 
$$ 
Throughout the proof we use the metric $\tau$ on $X$. We begin with the following

\begin{obs}
\label{o1}
Let $g_l\in G_{x_l}\setminus G_{e_{odd(l)}}$ and $z$ be a vertex of $g_l(F)$. 
Then there is a geodesic segment $\si=\ol{zz'}\subset (X,\tau)$ from $z$ to a point 
$z'\in \Ga_{odd}$, which intersects $T_{odd}$ orthogonally 
(at the point $z'$) and which is entirely contained in $F\cup g_l(F)$. 
For instance, if $g_l(F)\cap F=\{x_l\}$, 
then the segment $\si$  equals $\ol{zx_{l}}\cup \ol{x_l m_{odd(l)}}$. 

Moreover, unless $z\in g_l(F)\cap F$, the length of the segment $\si$ is 
strictly greater than $1$. In particular, $d(z, T_{odd})>1$.  
\end{obs}

Let $f:= h_j^{-1} g_i $ and assume that  $f(T_{odd})=T_{odd}$. Hence the edge $f(e_{odd(i)})\subset X$ intersects 
$T_{odd}$ orthogonally  in its midpoint $f(m_{odd(i)})$. The segment $f(e_{odd(i)})$ 
equals $\ol{zw}$ where $z:= h_j^{-1}(x_i)$ and $w:= f(x_{i+1})$ (if $i$ is odd),  
and $w:= f(x_{i-1})$ (if $i$ is even). In any case, $d(z, T_{odd})=1$. 

Suppose that $h_j\notin G_{odd(j)}$. 
Then, unless $z=x_j$, by applying Observation \ref{o1} to $l=j$, we obtain a contradiction 
with $d(z, T_{odd})>1$. Therefore $z=x_j$ and $i=j$, $f\in G_{x_i}$. Thus we apply 
Observation \ref{o1} to $l=j$, $g_l=f$ and the vertex $w\in f(F)$, and conclude that 
$w\in F$ as well. Then $f(e_{odd(i)})=e_{odd(i)}$ and hence 
$f=h_j^{-1} g_i \in \Phi_{odd}'$. 

Hence we conclude that $h_j\in G_{odd(j)}$. We now use the fact that 
$g_i^{-1}h_j(T_{odd})=T_{odd}$ and reverse the roles of $g_i$ and $h_j$. The same argument as above 
then implies that either 

\begin{itemize}
\item[(a)] $g_i\in G_{odd(i)}$, or  

\item[(b)] $i=j$, $f^{-1}\in G_{odd(j)}$.  
\end{itemize}

In Case (a), $h_j\in G_{odd(j)}$, $g_i\in G_{odd(i)}$ and thus $f\in \Phi_{odd}'$; 
in Case (b) $i=j, f \in G_{odd(j)}$ and thus $f\in \Phi_{odd}'$. \qed 

\begin{cor}
For each $g\in \Phi_{even}, h\in \Phi_{odd}$ we have $gT_{even}\cap T_{even}=\emptyset$, 
$hT_{odd}\cap T_{odd}=\emptyset$. 
\end{cor}

Recall the following definition:

\begin{definition}
A group $G$ is said to satisfy LERF property with respect to a subgroup $H\subset G$ if one of the 
following equivalent conditions holds:

(a)\qua For each finite subset $F\subset G\setminus H$ there exists a homomorphism $\phi\co  G\to \bar{G}$, 
where $\bar{G}$ is a finite group and $\phi(H)\cap \phi(F)=\emptyset$.  

(b)\qua For each finite subset $F\subset G\setminus H$ there exists a finite index subgroup $G'\subset G$ so 
that $H\subset G'$ and $F\cap G'=\emptyset$. 
\end{definition}

Note that for $H=\{1\}$ the above definition amounts to residual finiteness of $G$.

\begin{cor}
There is an epimorphism $\phi\co  G\to \bar{G}$ where $\bar{G}$ is a finite group and $\phi(\Phi_{even})\cap 
\phi(G_{even})=\emptyset$,  $\phi(\Phi_{odd})\cap \phi(G_{odd})=\emptyset$. 
\end{cor}
\proof According to \cite{Wise}, the group $G$ satisfies the LERF property with respect to each 
quasi-convex subgroup. Thus there are finite quotients
$$
\phi'\co  G\to \bar{G}', \phi''\co  G\to \bar{G}''
$$
so that $\phi'(G_{even})\cap \phi'(\Phi_{even})=\emptyset, 
\phi'(G_{odd})\cap \phi'(\Phi_{odd})=\emptyset$. Then define the 
homomorphism $\phi=(\phi', \phi'')\co  G\to \bar{G}'\times\bar{G}''$ and let the group $\bar{G}$ be  the 
image of $\phi$. \qed
\medskip  

Let $q\co  X\to \bar{X}:= X/Ker(\phi)$ denote the quotient map;  
the group $G$ acts on the compact complex $\bar{X}$ through the quotient group $\bar{G}$. 
We let $\bar{T}_{even},  \bar{T}_{odd}$ denote the projections of the trees $T_{even}, T_{odd}$ 
to the complex $\bar{X}$.

\begin{lem}
Suppose that $g_i\in G_{x_i}, h_j\in G_{x_j}$ and $g_i \bar{T}_{even}\cap h_j\bar{T}_{even}\ne \emptyset$ 
(resp.\ $g_i \bar{T}_{odd}\cap h_j\bar{T}_{odd}\ne \emptyset$). Then $h_j^{-1} g_i\in \Phi_{even}'$ 
(resp.\ $h_j^{-1} g_i\in \Phi_{odd}'$). 
\end{lem}
\proof If $g_i \bar{T}_{even}\cap h_j\bar{T}_{even}\ne \emptyset$ then 
$g_i \bar{T}_{even}=h_j\bar{T}_{even}$. It follows 
that $k h_j^{-1} g_i \in G_{even}$ for some $k\in Ker(\phi)$; thus $\phi(h_j^{-1} g_i)\in 
\phi(G_{even})$ which implies that $h_j^{-1} g_i\in  \Phi_{even}'$. 
The argument for $\bar{T}_{odd}$ is the same. 
\qed\medskip 

The graphs $\bar{T}_{even}, \bar{T}_{odd}$ determine finite subsets $S_{even}, S_{odd}$ 
of the set $Edges(\bar{X})$ consisting 
of those edges in $\bar{X}$ which intersect $\bar{T}_{even}, \bar{T}_{odd}$ nontrivially. 
Let $\xi, \eta$ denote the characteristic functions of the subsets $S_{even}, S_{odd}\subset 
Edges(\bar{X})$, normalized to have unit norm in the (finite-dimensional) Hilbert space 
$H:=L^2(Edges(\bar{X}))$. The group $G$ acts on $H$ by precomposition. 
We let $V\subset H$ denote the span of the subset $G\cdot \{\xi, \eta\}\subset H$ and let 
$p$ be the dimension of $V$. 

\begin{cor}
\label{ortho}
{\rm(1)}\qua The subgroups 
$G_{even}, G_{odd}$ fix the vectors $\xi, \eta\in V$ respectively.  

{\rm(2)}\qua The set 
$$
\Si=\{g^*(\xi), g^*(\eta): g\in G_{x_1}\cup \ldots G_{x_n}\}$$
is an orthonormal system in $V$. 

{\rm(3)}\qua For all $g, h\in G_{x_1}\cup \ldots G_{x_n}$: 

{\rm(a)}\qua $g^*(\xi)\ne h^*(\eta)$. 

{\rm(b)}\qua $g^*(\xi)=h^*(\xi)$ (resp.\ $g^*(\eta)=h^*(\eta)$) iff $h^{-1}g\in \Phi'_{even}$ 
(resp.\ $h^{-1}g\in \Phi'_{odd}$).  
\end{cor}

\noindent Let $Lk^\circ_{x_i}$ denote the vertex set of the link of $x_i$ in $X$. 
  
\begin{cor}
\label{CN}
The representation $G\acts V$ contains subrepresentations 
$$
G_{x_i}\acts Vect(Lk^\circ_{x_i}),$$
 so that 
the orthonormal vectors $G_{x_i}\cdot \{\xi, \eta\}$ are identified with the vectors   
 $G_{x_i}\cdot Lk^\circ_{x_i}$. 
\end{cor}
\proof Let us consider the case when $i$ is odd, since the other case is analogous. 
Observe that the stabilizer in $G_{x_i}$ of the vector $\xi\in V$ (resp.\ $\eta\in V$) 
is the group $G_{x_i}\cap G_{even}=G_{e_{even(i)}}$ (resp.\ 
$G_{x_i}\cap G_{odd}=G_{e_{odd(i)}}$); the stabilizer in $G_{x_i}$ of the vector $f_i^-$ is 
$G_{e_{even(i)}}$. Thus we construct an isometric embedding 
$Vect(Lk^\circ_{x_i})\to V$ by sending $f_i^-, f_i^+$ to $\xi, \eta$ respectively, 
and then extending this map equivariantly 
to the orthonormal basis $Lk^\circ_{x_i}$. \qed

\section{Hyperbolic trigonometry}
\label{trig}

\begin{figure}[ht!]\small\anchor{F0}
\psfrag {m}{$m$}
\psfrag {x}{$x$}
\psfraga <-2pt,0pt> {a}{$a_n$}
\psfrag {b}{$b_n$}
\psfrag {r}{$r_n$}
\psfrag {rho}{$\rho_n$}
\psfraga <-2pt,-2pt> {O}{$O$}
\psfrag {F}{$F$}
\psfraga <0pt,2pt> {p}{$\pi/n$}
\centerline{\includegraphics[width=2.5in]{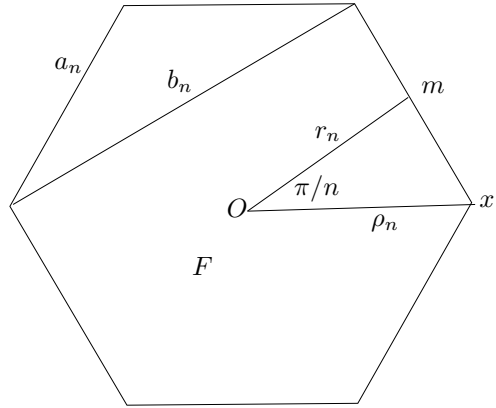}}
\caption{Geometry of a hyperbolic $n$--gon}
\label{F0}
\end{figure}

Consider a regular right-angled hyperbolic $n$--gon $F\subset \H^2$ ($n\ge 5$). 
Let $a_n$ denote its side-length, $\rho_n$ the radius of the circumscribed circle, 
$r_n$ the radius of the inscribed circle, $b_n$ the length of the shortest 
diagonal in $F$ (ie, a diagonal which cuts out a triangle from $F$);  
see \figref{F0}. We then have:
$$
\cosh(a_n)= 1+ 2 \cos(\frac{2\pi}{n}),  \quad 
\cosh(\frac{a_n}{2})=\sqrt{2}\cos(\frac{\pi}{n}),  \quad 
\cosh(b_n)= \cosh^2(a_n), $$
$$
\cosh(r_n)= \frac{1}{\sqrt{2}\sin(\frac{\pi}{n})}, \quad 
\cosh(\rho_n)= \cosh(r_n) \cosh(\frac{a_n}{2}). 
$$
\noindent Note that $a_n, b_n, r_n, \rho_n$ are strictly increasing functions of $n$. 

\begin{figure}[ht!]\small\anchor{F5}
\psfrag{x}{$x$}
\psfrag{y}{$y$}
\psfrag{zero angle}{zero angle}
\centerline{\includegraphics[width=2.5in]{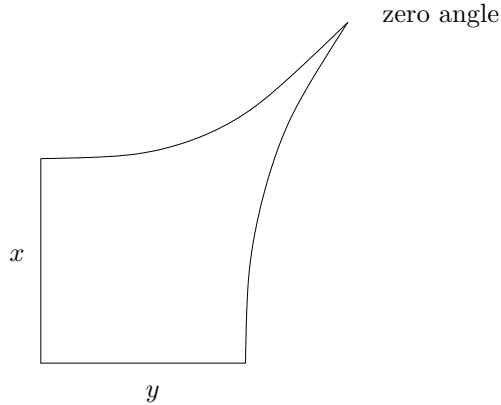}}
\caption{A Lambert quadrilateral}
\label{F5}
\end{figure}

Consider a {\em Lambert quadrilateral} $Q$ with one zero angle: $Q$ is a quadrilateral in $\H^2$ 
with one ideal vertex (where $Q$ has zero angle) and three finite vertices where the 
angles are $\pi/2$ (\figref{F5}). Let $x, y$ denote the lengths of the finite sides of $Q$. 
Then 
$$
\sinh(x)\sinh(y)=\cos(0)=1$$
 (see \cite[7.17.1]{Beardon}), or, equivalently
$$
\cosh^2(x)\cosh^2(y)= \cosh^2(x)+\cosh^2(y). 
$$ 
Thus, if we have two segments $E=\ol{xx'}, E'=\ol{x'x''}$ 
in $\H^2$ which intersect at the point $x'$ where they meet orthogonally,  
then the necessary and sufficient condition 
for $\Bis(E)\cap \Bis(E')=\emptyset$ is
$$
\cosh^2(|E|/2)\cosh^2(|E'|/2)\ge \cosh^2(|E|/2)+\cosh^2(|E'|/2),$$
equivalently
$$
\sinh(|E|/2)\sinh(|E'|/2)\ge 1. 
$$
We will refer to these inequalities as the {\em disjoint bisectors test}. 

\begin{lem}
\label{L1}
Suppose that $E, E'\subset \H^2$ are segments which meet orthogonally  
at a vertex, where $|E|=2\rho_n$ and $|E'|=a_n$. Then 
$$
\emptyset = \ol{\Bis(E)}\cap \ol{\Bis(E')}\subset \ol{\H}^2, 
$$
provided that $n\ge 7$; in case $n=6$ we have: 
$$
\ol{\Bis(E)}\cap \ol{\Bis(E')}\subset \geo \H^2, 
$$
is a point at infinity. 
\end{lem}
\proof Applying the disjoint bisectors test to $|E|=2\rho_6$, $|E'|=a_6$ 
we get the equality. Hence the bisectors meet at infinity in case $n=6$ and are 
within positive distance from each other if $n\ge 7$. \qed 

\begin{lem}
\label{L2}
Suppose that $E, E'\subset \H^2$ are segments which meet orthogonally  
at a vertex, where $|E|=b_n$ and $|E'|=a_n$, $n\ge 7$. Then 
$$
\emptyset = \ol{\Bis(E)}\cap \ol{\Bis(E')}\subset \ol{\H}^2. 
$$
\end{lem}
\proof Since $b_n\ge b_7, a_n\ge a_7$ it suffices to prove lemma in case $n=7$. Note that 
$a_7> a_6$ and
$$
b_7\approx 2.302366350 > 2\rho_6= 2.292431670
$$
Hence the assertion follows from Lemma \ref{L1}. \qed 

\medskip 
Below is another application of the {\em disjoint bisectors test.} 
Consider three segments $s, s', s''$ in $\H^3$ of the length $x, y, x$ respectively, 
which are mutually orthogonal and so that $s\cap s'=p, s'\cap s''=q$, $s'=\ol{pq}$, 
see \figref{F8}.

\begin{cor}
\label{C1}
If $x=b_n, y=a_n, n\ge 5$,  then $\ol{\Bis(s)}\cap  \ol{\Bis(s'')}=\emptyset$. 
\end{cor}
\proof It suffices to prove the corollary for $n=5$. We first compute the 
length $z=2t$ of the segment $s'''$ coplanar to $s$ and $s'$ such that 
$\Bis(s)=\Bis(s''')$. By considering the Lambert's quadrilateral with angle $\phi$ we get:
\begin{figure}[ht!]\small\anchor{F8}
\psfrag {p}{$p$}
\psfrag {x}{$x$}
\psfraga <3pt,-2pt> {s}{$s$}
\psfrag {h}{$h$}
\psfrag {y}{$y$}
\psfraga <-3pt,12pt> {z}{$z$}
\psfrag {x}{$x$}
\psfrag {s'}{$s'$}
\psfrag {s'''}{$s'''$}
\psfrag {t=z/2}{$t=z/2$}
\psfraga <-3pt,0pt> {q}{$q$}
\psfraga <-3pt,0pt> {f}{$\phi$}
\psfrag {s"}{$s''$}
\centerline{\includegraphics[width=2.5in]{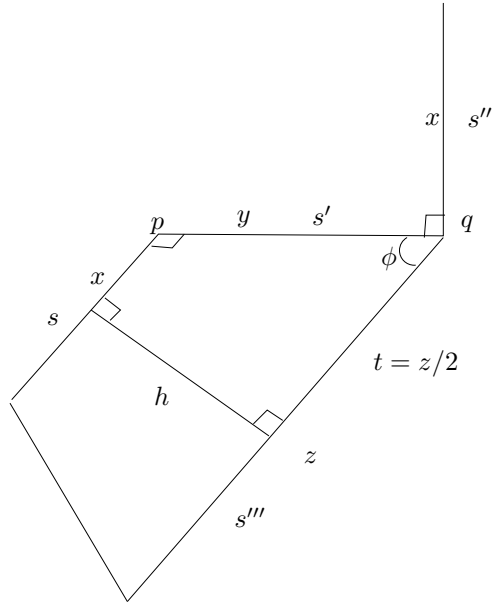}}
\caption{Three orthogonal segments}
\label{F8}
\end{figure}
$$
\left\{
\begin{array}{c}
\cosh(t)\sin(\phi)=\cosh(x/2)\\
\cosh(y)\sin(\phi)=\cosh(h)\\
\sinh(x/2)\sinh(h)=\cos(\phi)
\end{array}
\right.
$$ 
$$
\cosh(t) =\frac{\cosh(y)\cosh(x/2)}{\cosh(h)}.\leqno{\rm Thus}
$$
On the other hand, the last two equations in the above system imply that 
$$
\cosh(h)= \frac{\cosh(x/2)}{\sqrt{ \sinh^2(x/2) + 1/\cosh^2(y)}}. 
$$
Therefore
$$
\cosh(z/2)= \sqrt{ 1+ \sinh^2(x/2)\cosh^2(y)}. 
$$
By applying the {\em disjoint bisectors test} to $s''', s''$ we get:
$$
\Bis(s)\cap \Bis(s'')=\emptyset \iff \sinh^2(x/2)\cosh(y)\ge 1. 
$$

Lastly, we have: 
$$
\cosh(a_5)\approx 1.6180, \quad \sinh(b_5/2)\approx 1.85123.
$$
Therefore $\sinh^2(b_5/2)\cosh(a_5)>1$. \qed

\section{Quasi-isometric maps of polygonal complexes}

Suppose that $X$ is a simply-connected 2--dimensional regular cell complex which is 
equipped with a path-metric so that:

1. Each face is isometric to a right-angled regular $n$--gon in $\H^2$ (of course, $n\ge 5$). 

2. The complex $X$ is {\em negatively curved}, ie, for each vertex $x\in X$ 
the length of the shortest embedded loop in $Lk_{x}(X)$ is at least $2\pi$. 

\begin{thm}
\label{qi}
Suppose that $\mu\co  X\to \H^p$ is a continuous map which is a (totally-geodesic) isometric embedding 
on each face of $X$. We also assume that for each pair of faces $F', F''\subset X$ which intersect nontrivially 
a common face $F\subset X$, we have:
$$
\Span(\mu(F')) \perp \Span(\mu(F'')). 
$$
Then $\mu$ is a quasi-isometric embedding. 
\end{thm}
\proof Throughout the proof we will be using the notation $a_n, b_n, \rho_n$ for various 
distances in a regular right-angled hyperbolic $n$--gon, see Section \ref{trig}. 

Since the inclusion $X^{(1)}\embed X$ is a quasi-isometry, 
it suffices to check that $X^{(1)} \stackrel{\mu}{\to} \H^p$ is a 
quasi-isometric embedding. Since $\mu$ is 1--Lipschitz, it is enough to 
show that $d(\mu(z), \mu(w))\ge C\cdot d(z, w)$ for some $C=C(X)>0$ and all 
$z, w\in X^{(0)}$.  We first give a proof in case $n\ge 6$ and then explain how to modify it 
for $n=5$.

Let $\tilde\ga\subset X^{(1)}$ be an (oriented)  geodesic segment connecting 
$z$ to $w$. 
We start by replacing (in case when $n$ is even) each subsegment of $\t\ga$ 
connecting antipodal points in a face $F$ of $X$ with a geodesic segment 
within $F$. We will call the resulting (oriented) curve $\ga\subset X$. 
Clearly, 
$$
\Length(\mu(\ga))\le \Length(\mu(\t\ga)) \le \frac{n}{2}\Length(\mu(\ga)),
$$ 
so it suffices to get a lower bound on $\Length(\mu(\ga))$. 
We will refer to the edges of $\ga$ connecting antipodal points of faces as 
{\em diagonals} in $\ga$. 

\begin{figure}[ht!]\small\anchor{F1}
\psfrag {F}{$F$}
\psfrag {F'}{$F'$}
\psfrag {x}{$x$}
\psfrag {x'}{$x'$}
\psfrag {x"}{$x''$}
\psfrag {g}{$\gamma$}
\psfraga <-3pt,0pt> {gt}{$\tilde\gamma$}
\centerline{\includegraphics[width=3.5in]{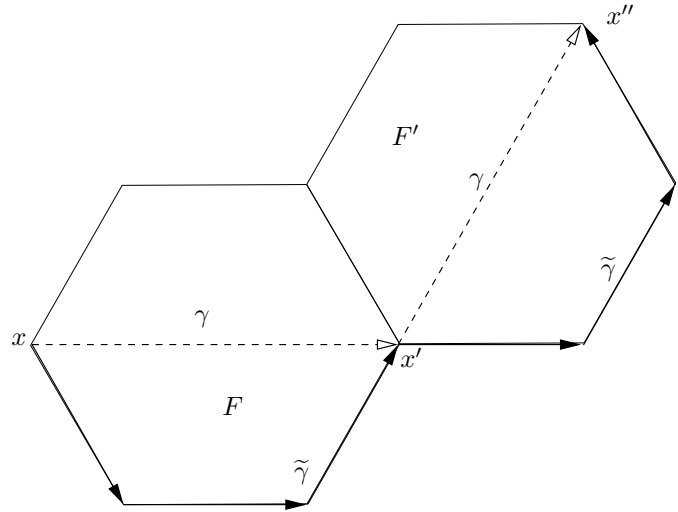}}
\caption{Diagonals}
\label{F1}
\end{figure}

\begin{rem}
Suppose that $\ol{x x'}, \ol{x' x''}$ are (distinct) diagonals in $\ga$, 
contained in faces $F, F'$ respectively. Then $F\cap F'= \{x'\}$: Otherwise 
 $\t\ga\subset X^{(1)}$ would  not be a geodesic as there exists a shorter path along 
 the boundaries of $F$ and $F'$; see \figref{F1}. 
In particular, $\mu(\ol{x x'}) \perp  \mu(\ol{x' x''})$ in $\H^p$. 
\end{rem}

We will regard $\ga$ as concatenation of consecutive segments $e_0, e_1,\ldots$. 
We define a collection $\b\e (\ga)$ of {\em bisected edges} $E_i$ in $\ga$ inductively as  follows:

(1)\qua Let $E_0=e_0\subset \ga$ be the first edge of $\ga$. 

(2)\qua Suppose that $E_i=e_j\subset \ga$ was chosen, $i\ge 0$. We will take as 
$E_{i+1}=e_k$, $k> i$, 
the first edge on $\ga$ following $E_i$  which satisfies two properties:

(a) If $e_i\cap e_k\ne \emptyset$ then either $e_i$ or $e_k$ is a diagonal. 

(b) $e_i, e_k$ do not belong to a common face in $X$. 

\medskip 

\begin{prop}
\label{P3}
Suppose that $n\ge 6$. Then the edges $E_i, i=0, 1, \ldots$ in $\ga$ satisfy 
the following:

{\rm(1)}\qua There exists a constant $c=c(X)$ such that $d(E_i, E_{i+1})\le c$ in $X$. 

{\rm(2)}\qua $\ol{\Bis(\mu(E_i))} \cap \ol{\Bis(\mu(E_{i+1}))}\subset \ol{\H}^p$ is 
empty unless $n=6$ and either $E_i$, or $E_{i+1}$ is not a diagonal. 
In case $n=6$ and at least one of these segments is not a diagonal, 
the bisectors $\Bis(\mu(E_i)), \Bis(\mu(E_{i+1})) \subset \H^p$ 
are disjoint but have a common ideal point in $\geo \H^p$. 

{\rm(3)}\qua For all edges $e_j\subset \ga$ between $E_i, E_{i+1}$, their images 
$\mu(e_j)$ are disjoint from  $\Bis(\mu(E_i))\cup \Bis(\mu(E_{i+1}))$. 
\end{prop}

\proof (1)\qua It is clear from the construction, that $E_i, E_{i+1}$ are separated by 
at most $n/2$ edges on $\ga$. Hence the first assertion follows. 

\begin{figure}[ht!]\small\anchor{F2}
\psfrag {E}{$E_i$}
\psfrag {E1}{$E_{i+1}$}
\psfrag {F}{$F$}
\psfrag {F'}{$F'$}
\psfrag {x}{$x$}
\psfrag {x'}{$x'$}
\psfrag {x"}{$x''$}
\psfrag {g}{$\gamma$}
\psfraga <-2pt,0pt> {gt}{$\tilde\gamma$}
\centerline{\includegraphics[width=3.5in]{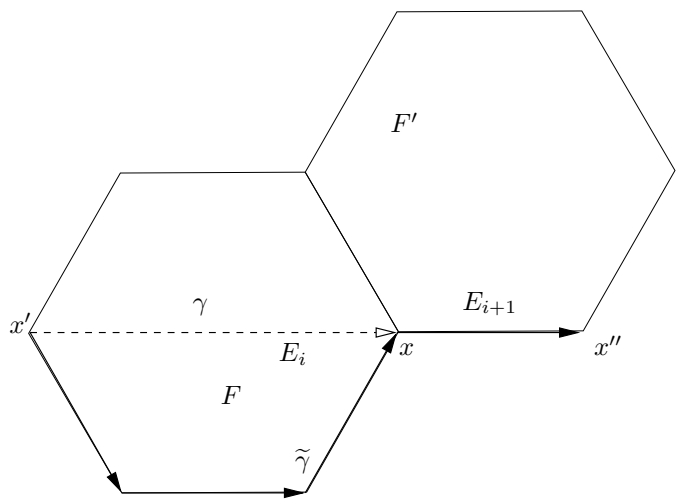}}
\caption{Diagonal and an edge}
\label{F2}
\end{figure}

(2)\qua There are several cases we have to consider. 
\medskip

(a)\qua Suppose that either $E_i$ or $E_{i+1}$ is a diagonal 
(see \figref{F2}) of the length $2\rho_n$ in the notation of section \ref{trig}. 
Then these segments share a common vertex $x'$ and it follows 
that $\mu(E_i) \perp  \mu(E_{i+1})$ (see the Remark above). The worst 
case occurs when $n=6$ and one of the segments is an edge of a face of $X$: 
The bisectors $\Bis(\mu(E_i))$,  $\Bis(\mu(E_{i+1}))$ 
are disjoint in $\H^p$ but have a common ideal point (see Lemma \ref{L1}). 
Since, as $n$ increases, both side-lengths and lengths of diagonals in regular 
right-angled $n$--gons in $\H^2$ strictly increase, it follows that 
$$
\ol{\Bis(\mu(E_i))} \cap \ol{\Bis(\mu(E_{i+1}))}=\emptyset, \forall n\ge 7. 
$$

(b)\qua Consider now the case when neither $E_i=e_k$ nor $E_{i+1}$ is a diagonal, 
$E_i$ is contained in a face $F$ and there exists at least one edge (say, 
$e_{k+1}$) between $E_i, E_{i+1}$ which is contained in the face $F$. Then, 
by the construction, $E_{i+1}=\ol{x' x''}$ is not contained in $F$ but shares  
the common point $x'$ with $F$. Thus $\mu(E_{e+1})\perp \mu(F)$. 

\begin{figure}[ht!]\small\anchor{F3}
\psfraga <0pt,-4pt> {E}{$E_i$}
\psfrag {Bis}{$\Bis(E_i)$}
\psfrag {E1}{$E_{i+1}$}
\psfraga <-4pt,-4pt>  {F}{$F$}
\psfrag {F'}{$F'$}
\psfraga <-2pt,0pt>  {y}{$y$}
\psfrag {x}{$x$}
\psfrag {x'}{$x'$}
\psfraga <-2pt,0pt>  {x"}{$x''$}
\psfrag {g}{$\gamma$}
\psfraga <-2pt,0pt> {gt}{$\tilde\gamma$}
\centerline{\includegraphics[width=3.5in]{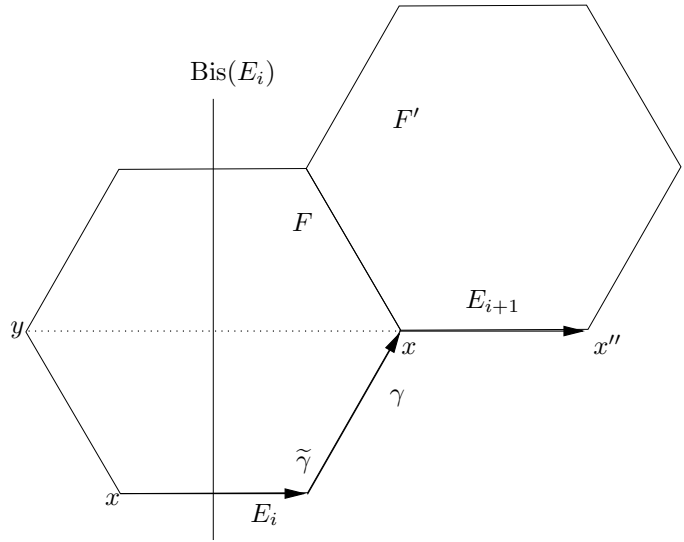}}
\caption{Bisectors}
\label{F3}
\end{figure}

Observe now that {\em there is a vertex $y\in F$ such that the segments $\ol{yx'}$ and 
$E_i$ have the same bisector in $F$.} To find this vertex simply apply the reflection 
in $\Bis(E_i)$ to the vertex $x'$: This symmetry preserves 
$F$ and sends the vertex $x'$ to a vertex $y\in F$.  See \figref{F3}. Since $x, x'\in F$ 
are not antipodal, $\ol{yx'}$ {\em is not an edge of $F$}. 

Clearly, $\Bis(\mu(\ol{yx'}))= \Bis(\mu(E_i))$. Hence the problem reduces to verifying 
that the bisectors $\Bis(\mu(\ol{yx'}))$, $\Bis(\mu(E_{i+1}))$ 
(or their closures in $\ol{\H}^p$) are disjoint. We note 
that in case $n=6$ the segment $\ol{yx'}$ connects antipodal points 
in $F$; hence the proof in this case reduces to (a). 

Assume now that $n\ge 7$, then, since $\ol{yx'}$ is not an edge of $F$, $|xy'|\ge b_n\ge b_7$, 
$|x'x''|=a_n\ge a_7$ and  
$$
\ol{\Bis(\mu(E_i))} \cap \ol{\Bis(\mu(E_{i+1}))}=\emptyset, 
$$
follows from Lemma \ref{L2}. 

(c)\qua The last case to consider is when $E_i, E_{i+1}$ are not diagonals and 
they are separated by exactly one edge $e\subset \ga$ 
(this edge cannot be a diagonal in this case), which is not contained 
in a common face with $E_i$ nor with $E_{i+1}$.  Then the edges 
$$
\mu(E_i), \mu(e), \mu(E_{i+1})\subset \H^p
$$
intersect orthogonally. The lengths of these edges are equal to $a(n)\ge a(6)$. Hence, 
(as in Case (b)) we replace $\mu(E_i)$  with a segment $s$ of the length 
$\ge 2\rho_6$ which meets $\mu(E_{i+1})$ orthogonally at the point 
$\mu(e)\cap \mu(E_{i+1})$. Therefore, by applying again 
Lemma \ref{L1}, the bisectors $\Bis(\mu(E_i)), \Bis(\mu(E_{i+1}))$ are 
disjoint; their closures in $\ol{\H}^p$ are disjoint provided $n\ge 7$. 

This proves the second assertion of the Proposition. The third assertion is 
clear from the construction: For instance, in Case (b) the edges $e_j$ 
between $E_i, E_{i+1}$ are all contained in the face $F$. 
Therefore they are disjoint from the bisector of $E_i$ within $F$, 
which implies the assertion about their images in $\H^p$. 
On the other hand, the edge $\mu E_{i+1}$ is orthogonal to $\mu(F)$, hence 
$\Bis(\mu E_{i+1})$ is disjoint from $\mu F$. \qed   

\medskip
Now, let us finish the proof that $\mu$ is a quasi-isometric embedding. Suppose that 
$\ga\subset X$ has length $L$, then the subset $\b\e(\ga)$ consists of 
$\ell\approx L/c$ {\em bisected edges} $E_i$ (here $c$ is the constant 
from Proposition \ref{P3}, Part 1). Hence the geodesic segment 
$\ga^*=\ol{\mu(z)\mu(w)}$ in $\H^p$ connecting the end-points 
of $\mu(\ga)$ crosses $\ell$ bisectors 
$\Bis(\mu E_i)$. In case $n\ge 7$, the bisectors   $\Bis(\mu E_i)$, 
$\Bis(\mu E_{i+1})$ are separated by distance $\delta=\delta(X)>0$, hence 
the length of $\ga^*$ is at least $\ell\delta$. Since $\ell\approx L/c$, 
we conclude that $d_{\H^p}(x, y)\ge Const\cdot L/\delta$. 
It follows that $\mu$ is a quasi-isometry. 

Now, consider the exceptional case $n=6$. We claim that for each $i$ 
the intersection points $\ol{\Bis(\mu E_i)}\cap \ol{\Bis(\mu E_{i+1})}$ 
and  $\ol{\Bis(\mu E_{i+2})}\cap \ol{\Bis(\mu E_{i+1})}$ are distinct. 
Given this, instead of the collection $\b\e(\ga)$ we would consider 
the collection of edges $E_i\in \b\e(\ga)$ for {\em even $i$}, then 
$\ol{\Bis(\mu E_i)}\cap \ol{\Bis(\mu E_{i+2})}=\emptyset$ for 
all even $i$ and we are done by the same argument as for $n\ge 7$.  

\medskip 
{\bf Case I}\qua We begin with the case when $E_i\subset F_i, 
E_{i+2}\subset F_{i+2}$ are diagonals and $E_{i+1}\subset F_i$ is not. 
(Here $F_i$ are faces of $X$.) 
Then $E_i\cap E_{i+1}$,  $E_{i+2}\cap E_{i+1}$ 
are the end-points of $E_{i+1}$. Therefore 
$$
\xi_i=\ol{\Bis(\mu E_i)}\cap \ol{\Bis(\mu E_{i+1})}\in \geo \Span(\mu(E_i)\cup  
\mu(E_{i+1})), 
$$
$$
\xi_{i+1}=\ol{\Bis(\mu E_{i+1})}\cap \ol{\Bis(\mu E_{i+2})}\in 
\geo \Span(\mu(E_{i+1})\cup \mu(E_{i+2})).  
$$
However, by the assumptions on $\mu$, 
$$
\Span(\mu(F_{i+2}))\perp \Span( \mu(F_i))
$$ 
Thus 
$$
 \Span(\mu(E_i)\cup \mu(E_{i+1}))\cap \Span(\mu(E_{i+1})\cup \mu(E_{i+2})) 
= \Span(\mu(E_{i+1})).  
$$
Since it is clear that $\xi_i\notin \geo  \Span(\mu(E_{i+1}))$, we conclude 
that $\xi_i\ne \xi_{i+1}$ and the assertion follows. 

\begin{figure}[ht!]\small\anchor{F4}
\psfrag {E}{$E_{i}$}
\psfrag {E1}{$E_{i+1}$}
\psfraga <-2pt,2pt> {E2}{$E_{i+2}$}
\psfrag {F}{$F_{i}$}
\psfrag {F1}{$F_{i+1}$}
\psfrag {F2}{$F_{i+2}$}
\psfrag {D}{$D_{i}$}
\psfrag {D2}{$D_{i+2}$}
\centerline{\includegraphics[width=3in]{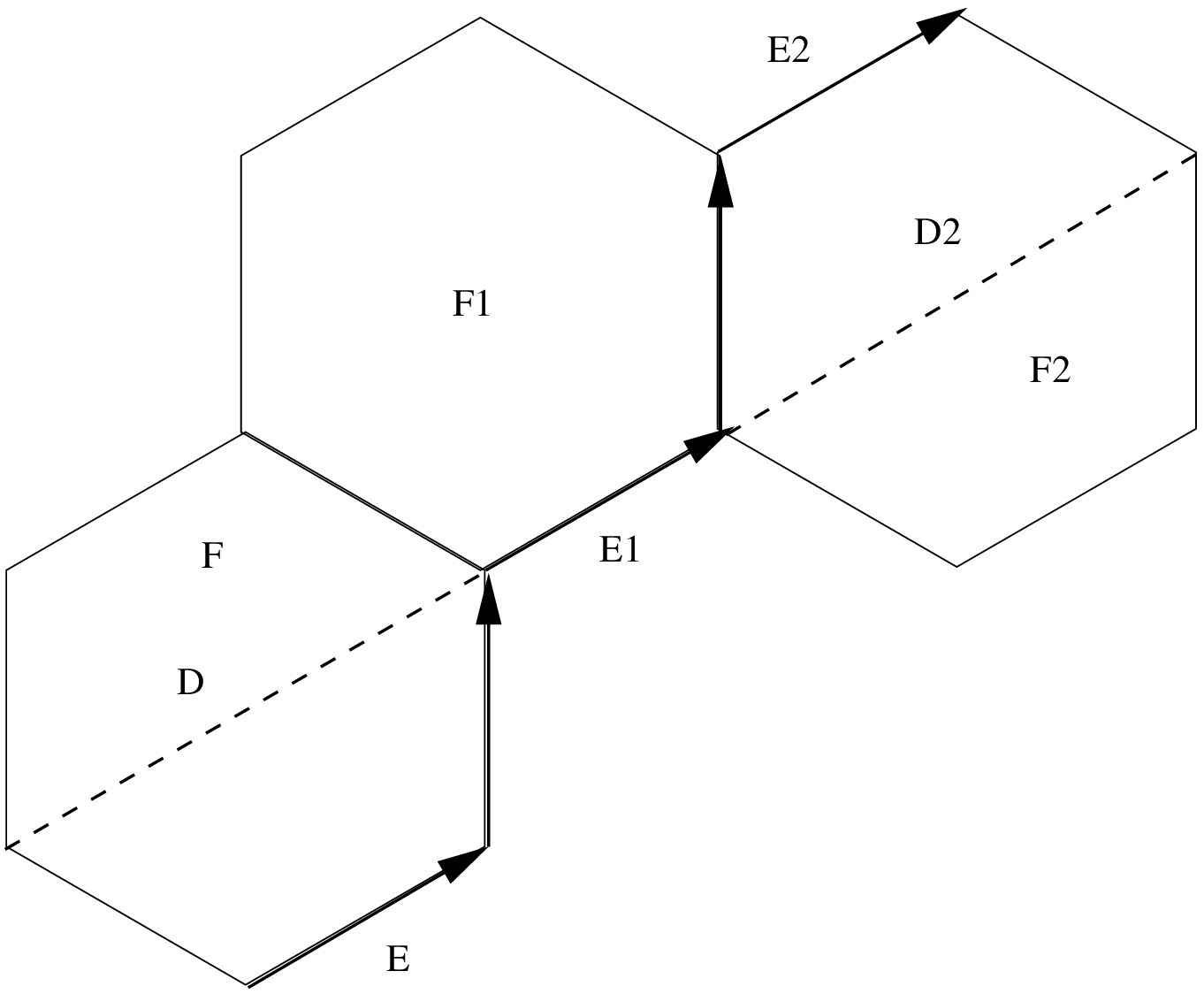}}
\nocolon\caption{}
\label{F4}
\end{figure}

\medskip 
We will reduce the case of a general triple of edges $E_{i}, E_{i+1}, E_{i+2}$ 
to the Case I discussed above. We consider only one other case, the 
arguments in the rest of cases are identical to: 

\medskip   
{\bf Case II}\qua Suppose that  pairs of edges $E_{i}\subset F_i, E_{i+1}\subset F_{i+1}$, 
$E_{i+2}\subset F_{i+2}$ are as in \figref{F4}. We take the diagonals 
$D_i\subset F_i, D_{i+2}\subset F_{i+2}$ so that
$$
\Bis(\mu D_i)= \Bis(\mu E_i), \Bis(\mu D_{i+2})= \Bis(\mu E_{i+2}). 
$$
Now the proof reduces to the Case I. 

\begin{figure}[ht!]\small\anchor{F6}
\psfrag {E1}{$E_{1}$}
\psfrag {E2}{$E_{2}$}
\psfrag {E3}{$E_{3}$}
\psfrag {F1}{$F_{1}$}
\psfrag {F2}{$F_{2}$}
\psfrag {F3}{$F_{3}$}
\psfrag {D1}{$D_{1}$}
\psfraga <-2pt,0pt> {D3}{$D_{3}$}
\psfrag {e}{$e$}
\psfrag {e'}{$e'$}
\centerline{\includegraphics[width=2truein]{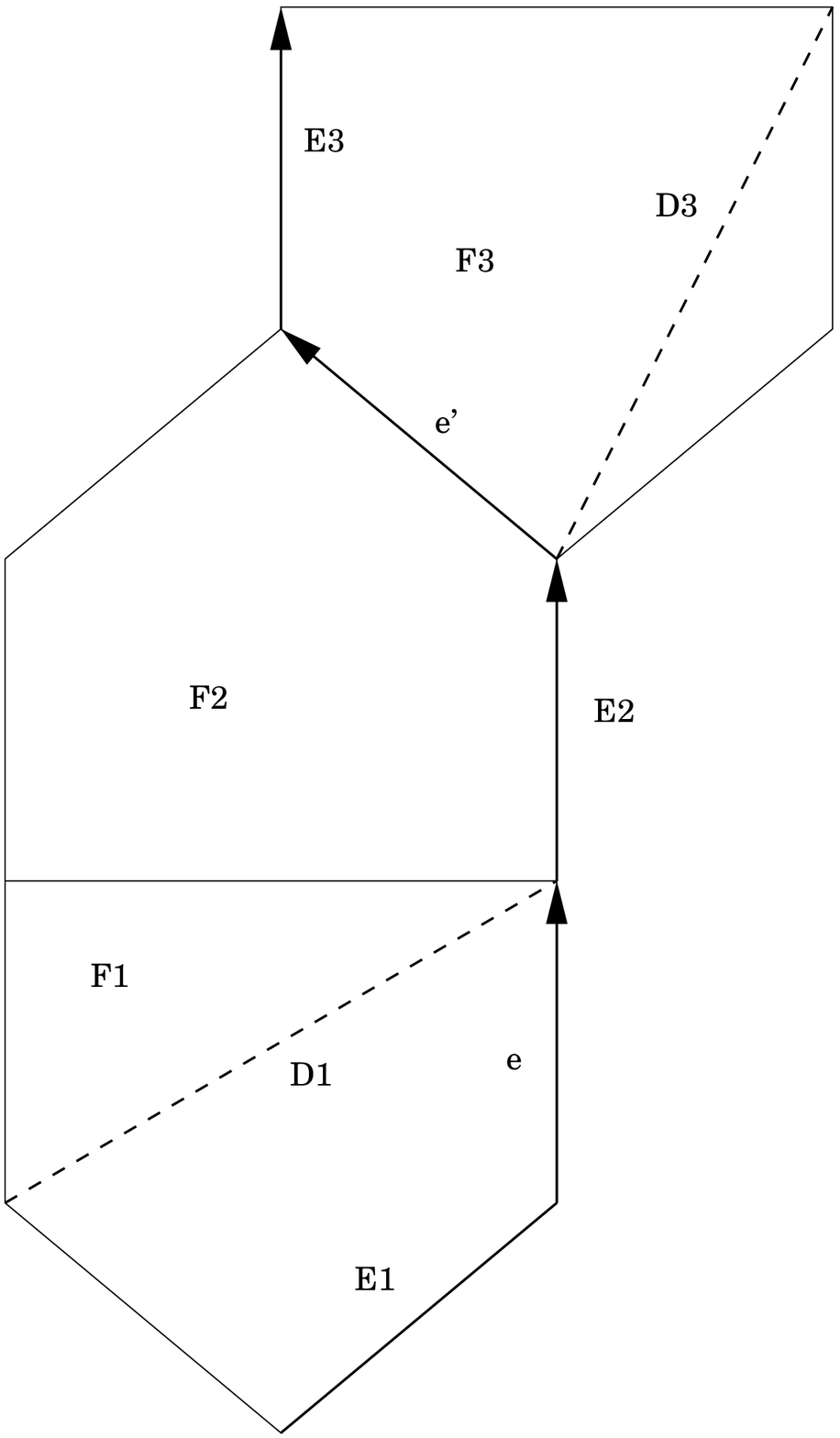}}
\nocolon\caption{}\label{F6}
\end{figure}

\medskip
Finally, consider the case of pentagons (ie $n=5$). We define 
the collection $\b\e(\ga)$ of {\em bisected edges} $E_i$ as before. 
Let $E_i, E_2, E_3$ be consecutive bisected edges. 
We will see that $\Bis(\mu E_1)\cap \Bis(\mu E_3)= \emptyset$. Since $n=5$  
we necessarily have: $E_2$ is separated by a unique edge $e\subset \ga$ 
from $E_1$ and by a unique edge $e'\subset \ga$ from $E_3$, see \figref{F6}. 

Note that it could happen that 
there is no face $F_1$ which contains $E_1, e$, nor  a face $F_3$ which contains $E_3, e'$. 
However, in $\H^p$ there exists a unique regular right angled pentagon which contains 
the edges $\mu(E_1), \mu(e)$ (resp.\ $\mu(E_3), \mu(e')$) in its boundary. Hence we will give 
a proof pretending that the corresponding face already exits in $X$. 
Observe that, similarly to our discussion above, the diagonals 
$D_1\subset F_1$, $D_3\subset F_3$ have the property that $\Bis(\mu D_i)= \Bis(\mu E_i)$, 
$i=1, 3$. 
Thus it suffices to consider the triple of pairwise orthogonal segments: 
$\mu(D_1), \mu(E_2)$ and  $\mu D_3$ in $\H^p$. The length of $\mu(E_2)$ 
equals $a_5$, the lengths of $D_1, D_3$ are equal to $b_5$, hence
$$
\ol{\Bis(\mu D_1)}\cap \ol{\Bis(\mu D_3)}=\emptyset
$$
by Corollary \ref{C1}. \qed

\section{Proof of the main theorem}

In this section establish 

\begin{thm}
\label{main1}
Suppose that $n=2k$ is even. Then the group $G$ admits an isometric properly discontinuous  
convex-cocompact action $\rho\co  G\acts\H^p$, 
where $p<\infty$ depends on the polygon $\p$.  
\end{thm}
\proof 
Let $X$ be the universal cover of the polygon of groups $\p$. 
We first construct a representation $\rho\co  G\to \isom(\H^p)$ for a certain $p$. 
We then produce a quasi-isometric $\rho$--equivariant embedding $\mu\co  X\to \H^p$. 
From this, via Lemma \ref{qconvex}, it will follow that $\rho\co  G\acts \H^p$ is 
an  isometric properly discontinuous convex-cocompact action.

\medskip 
Let $Lk_{x}$ denote the link (in $X$) of the vertex $x$, similarly, 
let $Lk_{e}$ denote the link of the edge $e$. Recall that $Lk^{\circ}$ denotes  
the vertex set of a link $Lk$. The set $Lk_{x_i}^\circ$ contains 
two distinguished elements: $f_i^+, f_i^-$ which correspond 
to the directions from $x_i$ toward $x_{i+1}$ and toward $x_{i-1}$ 
respectively. We define the subsets  $Lk_{x_i}^\bul := Lk_{x_i}^\circ \setminus 
\{ f_i^+, f_i^-\}$.  

Observe that the directions $\eta_i\in Lk_{x_i}^\circ, 
\eta_{i+1}\in Lk_{x_{i+1}}^\circ$ belong to the boundary of a common face in $X$ 
if and only if there exists $g\in G_{e_i}$ so that
$$
\eta_i= g(f_i^\pm), \eta_{i+1}= g(f_{i+1}^\pm). 
$$

{\bf Step 1: Construction of $\rho$} 

\medskip 
It is clear that to construct a representation $\rho\co  G\to \isom(\H^p)$ 
(for some $p$) we have to produce a collection of faithful representations
$$
\rho_i\co  G_{x_i} \to \isom(\H^p)
$$
so that the following diagram commutes: 
$$
\begin{array}{ccccc}
G_{x_i} & \longleftarrow & G_{e_i} & \longrightarrow & G_{x_{i+1}}\\
~&  \searrow\hbox{\scriptsize$\rho_i$} & ~~~ & \hbox{\scriptsize$\rho_{i+1}$} \swarrow & ~ \\
~& ~& \isom(\H^p) & ~& ~
\end{array}
$$ 
Embed $F$ isometrically (as a convex, regular, right-angled polygon) in the hyperbolic 
plane $\H^2$. Via this embedding we will identify the directions 
$f_i^{\pm}\in Lk_{x_i}$ with the unit vectors 
$\ov{f_i^{\pm}}\in T_{x_i}\H^2$ which are tangent to the sides of $F$.  

In what follows we will adopt the following convention. 
Given a number $p$ and a totally-geodesic embedding $\H^2\subset \H^p$, we  
observe that the normal bundle $N(\H^2)$ of $\H^2$ in 
$\H^p$ admits a canonical flat orthogonal connection 
(invariant under the stabilizer of $\H^2$ in $\isom(\H^p)$). 
Thus, given normal vectors $\nu\in N_x(\H^2), \nu'\in N_{x'}(\H^2)$, 
we have well-defined scalar product $\nu \cdot \nu'$ 
and hence the notion of orthogonality $\nu \perp \nu'$. Set  
$V_i:= T_{x_i}\H^p$ and let $N_i\subset V_i$ denote the orthogonal complement to 
$T_{x_i}\H^2$. We define $R_i\in \isom(\H^p)$ to be the isometric reflection 
in the bisector of the edge $e_i$ of $F\subset \H^p$. Set 
$J_i:= R_{i-1}\circ \ldots \circ R_{1}$, 
for $i=2,\ldots,n+1$; observe that $J_{n+1}=Id$.

\begin{rem}
The fact that the identity $J_{n+1}=Id$ fails if $n$ is odd is one of the reasons why our 
construction requires $n$ to be an even number. An attempt to apply the constructions below to 
odd $n$ lead to a representation of a certain extension of the group $G$ rather than of $G$ itself. 
\end{rem}

\begin{prop}
\label{P2}
Suppose that $n=2k$ is even. 
Then there exists a natural number $p$ and a collection of faithful isometric 
linear actions $d\rho_i\co  G_{x_i} \acts V_i$, $i=1,\ldots,n$, so that the following hold:

{\rm(1)}\qua Each representation $d\rho_i$ contains a subrepresentaion $G_{x_i}\acts 
Vect(Lk^\circ_{x_i})$, so that the unit basis vectors 
$f^{\pm}_i\in Lk^\circ_{x_i}\subset Vect(Lk^\circ_{x_i})$ are identified with the vectors 
$\ov{f^{\pm}_i}\in V_i$.  

{\rm(2)}\qua Each reflection $R_i\co  \H^p \to \H^p$ induces an isomorphism of 
$\R G_{e_i}$--modules 
$$(V_i, d\rho_i(G_{e_i}))\to (V_{i+1}, d\rho_{i+1}(G_{e_i})).$$ 

{\rm(3)}\qua ``Orthogonality'': The spaces $Vect(Lk_{x_i}^\bullet)\subset N_i, 
Vect(Lk_{x_j}^\bullet)\subset N_j$ 
are mutually orthogonal, $|i- j|\ge 2, i, j\in \{1,\ldots,n\}$.  
If $j=i+1$ then we require orthogonality of 
the subspaces 
$$
Vect(Lk_{x_i}^\circ \setminus G_{e_i}\cdot \{ f_i^+, f_i^-\})\subset N_i, 
Vect(Lk_{x_{i+1}}^\circ \setminus G_{e_i}\cdot \{ f_{i+1}^+, f_{i+1}^-\}) 
\subset N_{i+1}.
$$ 
\end{prop}

\begin{rem}
(1)\qua The assumption that the number of sides of $F$ is even is used only in 
this part of the proof of the main theorem and very likely is just a technicality. 

(2)\qua The ``orthogonality'' property will be used to prove that the action $G\acts \H^p$ 
that we are about to construct,  is discrete, faithful and convex-cocompact. 
%On the other hand, a slight modification of this part of the 
% construction also works for $n=4$ provided that $G$ is virtually torsion-free (where 
% we consider a regular quadrilateral $F\subset \H^2$ with acute angles). 
\end{rem}

Before beginning the proof of the proposition we first make some observations 
(where we ignore the {\em orthogonality} issue). Suppose that we have constructed 
representations $d\rho_i$. We then ``fold'' these representations into a single orthogonal 
representation $G\acts V_1$ by composing 
each $d\rho_i$ with the composition of reflections $(R_{1})_*\circ \ldots \circ (R_{i-1})_*$, where 
$(R_j)_*$ is the isomorphism $O(V_{j+1})\to O(V_j)$ which is induced by $dR_{j}|_{x_{j+1}}$. 
Note that under the action $G\acts V_1$ the vectors $\ov f_1^{+}, \ov f_1^{-}$ are fixed by  
the ``odd'' and ``even'' subgroups $G_{odd}, G_{even}$, respectively. Moreover, the representation 
$G\acts V_1$ contains subrepresentations $G_{x_i}\acts Vect(Lk_{x_i}^\circ)$. 

Recall that in Corollary \ref{ortho} we have constructed a finite-dimensional orthogonal 
representation $G\acts V$ which satisfies the same properties as above: It contains unit vectors 
 $\xi, \eta$ fixed by $G_{even}, G_{odd}$ respectively, and it contains 
subrepresentations $G_{x_i}\acts Vect(Lk_{x_i}^\circ)$. Therefore, to construct the 
representations $d\rho_i$ we begin with the action $G\acts V$ (which we identify with an 
action $G\acts V_1$) and then ``unfold'' it (using compositions of reflections $R_i$) 
to a collection of representations $d\rho_i$. This is the idea of the proof of Proposition \ref{P2}. 

The reader familiar with {\em bending} deformations of representations of groups into 
$\isom(\H^p)$ will notice  that the ``folding'' and ``unfolding'' of representations 
discussed above is nothing but {\em $\pi$--bending}.  

\proof We let $G\acts V$ denote the orthogonal representation constructed in Corollary \ref{ortho}. 
According to Corollary \ref{CN}, the representation $G\acts V$ contains subrepresentations 
$G_{x_i}\acts Vect(Lk^\circ_{x_i})$. Let $p$ denote the dimension of $V$. 

Our goal is to construct isometries $\phi_i\co  V\to V_i:= T_{x_i} \H^p$. 
The actions $d\rho_i$ will be obtained by the conjugation:  
$$
d\rho_i\co  G_{x_i}\acts V_i := 
\phi_i \circ (G_{x_i}\acts V) \circ \phi_i^{-1}.$$ 

First take an arbitrary isometry $\phi_1\co  V\to V_1$ sending the unit vectors $\xi, \eta\in V$ to 
the vectors  $\ov f_1^-, \ov f_1^+$ respectively. Now define isometries $\phi_j$, $j=2,\ldots,n+1$ 
by 
$$
\phi_i= J_i\circ \phi_1,
$$ 
$$
\phi_{i+1}= R_i\circ \phi_i.\leqno{\hbox{ie}}
$$
Note that $\phi_{n+1}=\phi_1$. Define the action $d\rho_i\co  G_{x_i} \acts V_i$ by conjugating 
via $\phi_i$ the action $G_{x_i} \acts V$.

The group $G_{e_1}$ fixes the vector $\eta\in V$, 
hence $G_{e_1}$ also fixes the vector $\phi_1(\eta)= \ov f_1^+$. Thus
$$
dR_1\circ d\rho_1\restr_{G_{e_1}}=d\rho_2\restr_{G_{e_1}}.  
$$
The same argument shows that 
$$
dR_i\circ d\rho_i\restr_{G_{e_i}}=d\rho_{i+1}\restr_{G_{e_i}}.  
$$
for all $i$. This proves (1) and (2). In what follows we will identify the spaces 
$Vect(Lk^\circ_{x_i})$ with their images in $V_i$, $i=1,\ldots,n$. 

\medskip  
We will check that the sets 
$$
Lk_{x_i}^\circ \setminus G_{e_i}\cdot \{ f_i^+, f_i^-\}\subset N_i, 
Lk_{x_{i+1}}^\circ \setminus G_{e_i}\cdot \{ f_{i+1}^+, f_{i+1}^-\} 
\subset N_{i+1}
$$ 
are orthogonal to each other and will leave the remaining orthogonality assertion to the reader. 
Let 
$$
v\in Lk_{x_i}^\circ \setminus G_{e_i}\cdot \{ f_i^+, f_i^-\}\subset N_i, 
w\in Lk_{x_{i+1}}^\circ \setminus G_{e_i}\cdot \{ f_{i+1}^+, f_{i+1}^-\}\subset N_{i+1}.$$
In order to show that $v\perp w$ it suffices to verify that the corresponding vectors
$$
v, w\in \Si\subset V
$$
are distinct (recall that $\Si$ is an orthonormal system in $V$). If, say, 
$v\in G\cdot \xi, w\in G\cdot \eta$ then $v\ne w$. Hence we will consider the case 
$$
v=g^*\eta\in G_{x_i}\cdot \eta, w=h^*\eta\in G_{x_{i+1}}\cdot \eta .
$$ 
According to Corollary \ref{ortho}, if $g^*(\eta)= h^*(\eta)$ then $h^{-1}g\in \Phi'_{odd}$. 
In our case, $g\in G_{x_i}, h\in G_{x_{i+1}}$. Then $h^{-1}g\in \Phi'_{odd}$ means that either

(a)\qua $g, h$ do not have a common fixed vertex of $F$ and $g\in G_{odd(i)}, h\in G_{odd(i+1)}$, 

\no or 

(b)\qua $g, h$ fix the same vertex of $F$ and $g^*(\eta)=h^*(\eta)$. 

In case (a), if $i$ is odd then $g, h\in G_{e_i}$, and hence 
$g^*\eta=g^*\eta$ corresponds to the vectors $f_i^+\in N_i, f_i^-\in N_{i+1}$. Therefore 
the equality $v=w$ implies that 
$$
\phi_i(v)\in G_{e_i}\cdot \ov f_i^+, \phi_{i+1}(w)\in G_{e_i}\cdot \ov f_{i+1}^-.$$ 
If $i$ is even then both $g\in G_{e_{i-1}}, h\in G_{e_{i+1}}$ fix the vector $\eta$. Therefore the 
equality $v=w$ implies that 
$$
\phi_i(v)=\ov f_i^-\in G_{e_i}\cdot \ov f_i^-, 
\phi_{i+1}(w)= \ov f_{i+1}^+\in G_{e_i}\cdot \ov f_{i+1}^+.$$
In case (b), we can assume that, say, $h\in G_{x_i}\cap G_{x_{i+1}}=G_{e_i}$. 
If $i$ is odd then $h^*\eta=\eta$ and therefore $g^*\eta=\eta$. 
Therefore the equality $v=w$ implies that
$$
\phi_i(v)=f_i^+\in G_{e_i}\cdot f_i^+, \quad \phi_{i+1}(w)= f_{i+1}^-\in G_{e_i}\cdot f_{i+1}^-.
$$
Lastly, assume that $i$ is even. Then the vector $v=g^*(\eta)= w=h^*(\eta)$ 
corresponds to the vectors
$$
\phi_i(w)=h(\ov f_i^-)\in G_{e_i}\cdot \ov f_i^-\subset N_i, \quad 
\phi_{i+1}(w)= h(\ov f_{i+1}^+)\in G_{e_i}\cdot \ov f_{i+1}^+\subset N_{i+1}. 
$$
This proves the {\em orthogonality} assertion.  \qed 

\medskip 
Now, once we have constructed linear orthogonal representations 
$d\rho_i\co  G_{x_i} \acts V_i$, we extend them (by exponentiation) to 
isometric actions $\rho_i \co   G_{x_i} \acts \H^p$, which fix the points 
$x_i$, $i\in 1,\ldots, n$. 
Observe that for each $i$ the group $\rho_i(G_{e_i})$, resp.\  
$\rho_i(G_{e_{i-1}})$, fixes the edge $e_i$, resp.\ $e_{i-1}$, of the polygon 
$F\subset \H^2$, since $d\rho_i(G_{e_i})$, $d\rho_i(G_{e_{i-1}})$, 
fix the vectors $\ov{f_i^+}, \ov{f_i^-}$ respectively. 
Hence the reflection $R_i$ commutes with the groups 
 $\rho_i(G_{e_i})$ and $\rho_{i+1}(G_{e_i})$. Thus the representations 
$\rho_i, \rho_{i+1}\co  G_{e_i} \to \isom(\H^p)$ are the same. Therefore  
the representations $\rho_i$ determine an isometric action $\rho\co  G\acts \H^p$.

\eject
{\bf Step 2: Discreteness of $\rho$} 
\medskip

We will construct a $\rho$--equivariant continuous mapping $\mu\co X\to
\H^p$ satisfying the assumptions of Theorem \ref{qi}.  Since such
$\mu$ is necessarily a quasi-isometric embedding (Theorem \ref{qi}),
by applying Lemma \ref{qconvex}, we will conclude that $\rho\co G\acts
\H^p$ is properly discontinuous and convex-cocompact.

Recall that we have identified the face $F$ of $X$ with a regular right-angled hyperbolic 
polygon $F$ in $\H^2$, this defines the (identity) embedding $\mu\co  F\to \H^2$. Now, 
for each $g\in G$ we set
$$
\mu\restr_{gF}:= \rho(g)\circ \mu\restr_{F}. 
$$ 
Let us check that this mapping is well-defined:

(1)\qua If $g\in G_F$ then, by construction of $\rho$, $\rho(g)$ fixes the polygon $\mu(F)$ pointwise, 
hence $\mu\circ g\restr_{F}=\rho(g)\circ \mu\restr_{F}$ for $g\in G_{F}$. 

(2)\qua If $g\in G_{e}$, where $e$ is an edge of $\p$, then, by the construction,  
$\rho(g)$ fixes the edge $\mu(e)$ pointwise, 
hence $\mu\circ g\restr_{e}=\rho(g)\circ \mu\restr_{e}$ for $g\in G_{e}$.

(3)\qua The same argument applies to $g\in G_x$, for the vertices $x\in F$. 

Hence $\mu\co  X\to \H^p$ is well-defined and thus it is a $\rho$--equivariant,  
continuous mapping which is an isometric totally-geodesic embedding on each face of $X$. 

Lastly, we check the {\em orthogonality} condition required by Theorem \ref{qi}. 
By equivariance, it is clear that we only need to verify orthogonality for the faces 
$F', F''\subset X$ which are adjacent to the face $F\subset X$. 
We will see that this orthogonality 
condition will follow from the Assertion 3 of Proposition \ref{P2}. 
There are several cases which may occur, we will check one of 
them and will leave the rest to the reader. 

Suppose that $F'=g_i(F), F''=g_j(F), g_i\in G_{x_i}, G_j\in G_{x_j}$ and 
$$
F'\cap F= \{x_i\}, \quad F''\cap F= \{x_j\},
$$
where $|i-j|\ge 2$. Then 
$$
T_{x_i} (\mu F')\subset N_i, \quad T_{x_j} (\mu F'')\subset N_j
$$
and the vectors 
$$
dg_{i}(\ov{f}_i^{\pm})\ne \ov{f}_i^{\pm}, \quad dg_{j}(\ov{f}_j^{\pm})\ne \ov{f}_j^{\pm}
$$
span $T_{x_i} (\mu F')$ and $T_{x_j} (\mu F'')$ respectively. 
According to the Assertion 3 of Proposition \ref{P2}, we have:
$$
dg_{i}(\ov{f}_i^{\pm})\perp dg_{j}(\ov{f}_j^{\pm}). 
$$
(Recall that here the  orthogonality is defined modulo the parallel translation along curves  
in $\Span(F)$.) Since both $\Span(\mu F'), \Span(\mu F'')$ intersect $\Span(F)$ orthogonally, 
the geodesic segment $\ol{x_i x_j}\subset \Span(F)$ is orthogonal to both 
$\Span(\mu F')$ and $\Span(\mu F'')$. Therefore $\mu(F')\perp \mu(F'')$. \qed

\section{The odd case}
\label{odd}

In this section we will construct examples of negatively curved right-angled  
polygons of groups and their actions on $\H^p$ in the case of the odd number 
of sides. 

\begin{figure}[ht!]\small\anchor{F11}
\psfrag {G1}{$G_{x_4}=\Gamma_1\times\Gamma_2$}
\psfrag {G2}{$G_{x_2}=\Gamma_2\times\Gamma_3$}
\psfrag {G3}{$G_{x_3}=\Gamma_3\times\Gamma_4$}
\psfrag {G4}{$G_{x_4}=\Gamma_4\times\Gamma_5$}
\psfrag {G5}{$G_{x_5}=\Gamma_5\times\Gamma_1$}
\psfrag {Ga1}{$\Gamma_1$}
\psfrag {Ga2}{$\Gamma_2$}
\psfrag {Ga3}{$\Gamma_3$}
\psfrag {Ga4}{$\Gamma_4$}
\psfrag {Ga5}{$\Gamma_5$}
\centerline{\includegraphics[width=4truein]{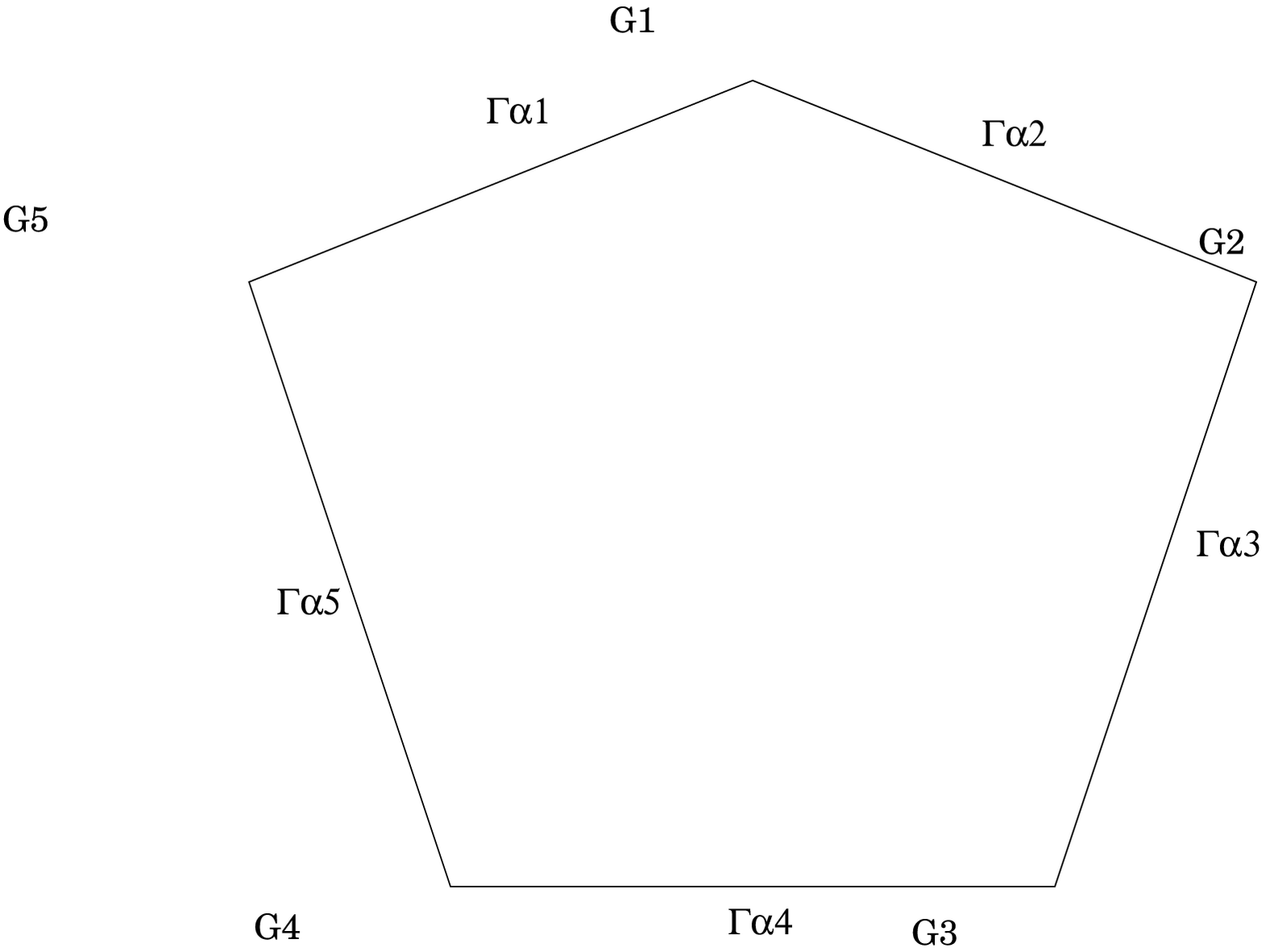}}
\nocolon\caption{}\label{F11}
\end{figure}

\medskip
We define the following polygon of groups. 
Suppose that we are given finite groups $\Ga_{1},\ldots,\Ga_{n}$.  
Let $F$ be an $n$--gon ($n\ge 3$). Below the indices $i$ are taken modulo $n$. 
We assign the group $G_{x_i}=\Ga_{i}\times \Ga_{{i+1}}$ 
to each vertex $x_i$ of $F$. We label each edge $e_i$ of $F$ by the group $\Ga_{i+1}$. 
The homomorphisms $\Ga_{i+1}\to   \Ga_{i}\times \Ga_{{i+1}}$,  
$\Ga_{i+1}\to   \Ga_{i+1}\times \Ga_{{i+2}}$ are the natural isomorphisms to the 
second and the first factor respectively. We set $G_F:=\{1\}$. In what follows, 
let $\p$ denote the resulting polygon of groups (see \figref{F11}) 
and set $G:= \pi_1(\p)$.

\begin{rem}
Note that $G$ is isomorphic to the {\em cyclic graph-product} of the groups $\Ga_i, i=1,\ldots,n$ 
(see \cite{Bridson-Haefliger} for detailed definition). Indeed, the group $G$ is generated 
by the elements of $\Ga_1,\ldots,\Ga_n$ subject to the relations:
$$
[g_i, h_{i+1}]=1, \forall g_i\in \Ga_i, h_{i+1}\in \Ga_{i+1}. 
$$ 
\end{rem}

The polygon of groups $\p$ is negatively curved provided that $n\ge 5$: For each 
vertex  $x=x_i\in \p$, the link of $x$ in the universal cover $X$ of $\p$ is the  
complete bipartite graph $K_{t_{i-1}, t_i}$, where $t_j:= |G_{e_{j}}|, j=1,\ldots,n$.

\begin{thm}
\label{graph}
Suppose that $n\ge 5$, and the polygon of groups $\p$ is as above. 
Then $G$ admits a properly discontinuous convex-cocompact action on 
$\H^p$ for certain $p=p(\p)$. 
\end{thm}

\proof We will assume that $n=2k-1$ is odd. Let $e_1,\ldots,e_n$ denote 
the edges of the polygon $\p$. 
First, we construct a homomorphism $\rho\co  G\to \isom(\H^p)$ for some $p$, 
which is faithful on each vertex group. 
Let $T$ denote the disjoint union  
$$
\bigcup_{i=1}^n Lk_{e_i},
$$
where the link of $e_i$ is taken in $X$. 
One can think of this set as the set of all flags: $(e, f)$, 
where $e$ is an edge in $F$ and $f$ is a face of $X$ 
containing $e$. Each group $G_{e_i}$ acts naturally on $Lk_{e_i}$ 
(since $G_{e_i}$ fixes the edge $e_i\subset X$). 
We extend this action to the {\em trivial action} on the rest of $T$. Thus we get an action 
$$
\prod_{i=1}^n G_{e_i} \acts T. 
$$
Observe that there is a tautological epimorphism
$$
G\to \prod_{i=1}^n G_{e_i}
$$
which sends each subgroup $G_{e_i}\subset G$ to the subgroup $G_{e_i}$ of the direct product. 
Hence $G$ acts on $T$ through the quotient group $\prod_{i=1}^n G_{e_i}$.

Let $W$ denote the Euclidean vector space $Vect(T)$ with 
the orthonormal basis $T$, and set $p:= dim(W)$, ie, 
$$
p=t_1+\ldots+t_n. 
$$ 
The set $T$ contains distinguished elements $f_1,\ldots,f_n$ 
consisting of the flags $(e_i, F)$. The dihedral group  
$D_n$ acts on $F$ and therefore on $\{f_1,\ldots,f_n\}$. 
We extend this action to the rest of $T$ (and hence to $W$) 
by the identity on $T\setminus \{f_1,\ldots,f_n\}$.  
The reflections $R_j\in D$ yield isometric involutions $I_j$ of $W$. 
Note that (since $n=2k-1$ is odd) the involution $R_{i+k}$ fixes 
$x_{i+1}$, hence $I_{i+k}$ permutes the vectors $f_{i}, f_{i+1}\in W$. 

Embed the polygon $F$ to $\H^2$ as a right-angled regular polygon 
$\mu(F)\subset \H^2$ as in the proof 
of Theorem \ref{main}. As before, we will identify $F$ and $\mu(F)$. 
We embed $\H^2$ into $\H^p$ as a totally-geodesic subspace. 
Let $m_i$ denote the midpoint of the edge $e_i\subset F$. 
The tangent space $W_i:=T_{m_i}\H^p$ contains 
a distinguished vector $\ov{f_i}$ which is the unit vector orthogonal to 
$e_i$ and directed inward $F$. 
Note that this vector is parallel (under the parallel transport along $e_i$) to the vectors 
$\ov{f}_i^+,  \ov{f}_{i+1}^-$ used in the proof of Theorem \ref{main}. 
Pick an arbitrary linear isometry $\psi_1\co  W\to W_1$ which 
sends $f_1$ to $\ov{f_1}$. 
% and $f_n$ to $\ov{f_1}^+$ (after the latter vector is translated into $W_1$ along $e_1$). 
By conjugating via $\psi_1$ we transport the linear action $G\acts W$ to a linear action 
$G\acts W_1$, exponentiating the latter action we get an isometric action $\rho_1\co  G\acts \H^p$.

We now proceed analogously to the proof of Theorem \ref{main}: 
Define linear maps $\psi_i\co  W\to W_i$ so that 
we have a commutative diagram:
$$
\begin{array}{ccccccc}
W & \stackrel{I_{k+1}}{\longrightarrow} & W & \stackrel{I_{k+2}}{\longrightarrow}&  
W ~~\ldots &\stackrel{I_{n+k+1}}{\longrightarrow}& W\\
~~\downarrow \hbox{\scriptsize$\psi_1$} & ~ & ~~\downarrow \hbox{\scriptsize$\psi_2$} & ~& \downarrow \hbox{\scriptsize$\psi_3$} & ~& 
~~\downarrow \hbox{\scriptsize$\psi_{n+1}$}\\
W_1 & \stackrel{R_{k+1}}{\longrightarrow} & W_2 & \stackrel{R_{k+2}}{\longrightarrow}&  
V_3 ~~\ldots &\stackrel{R_{n+k+1}}{\longrightarrow}& W_n\\
\end{array}
$$
\no Note that $\psi_i(f_i)= \ov{f_i}$ for all $i$. 
Indeed, $\psi_1(f_1)=\ov{f_1}$ by construction. Suppose that 
$\psi_i(f_i)=\ov{f_i}$. Then  
$$
\psi_{i+1}(f_{i+1})= R_j\circ \psi_i \circ I_j(f_{i+1})=  
R_j\circ \psi_i(f_i)= R_j(\ov{f_i})= \ov{f_{i+1}}
$$
where $j=i+k$. 

Observe that $\psi_{n+1}\ne \psi_1$. However, 
$I_{n+k+1}\circ \ldots \circ I_{k+1}= I_1$ commutes with $G_{e_1}\acts W$ and 
$R_{n+k+1}\circ \ldots \circ R_{k+1}= R_1$ commutes with $\rho_1(G_{e_1})\acts \H^p$. Hence 
$$
\rho_1=\rho_{n+1}\co  G_{e_1}\to \isom(\H^p).$$
 It remains to verify that for each $i$ the groups
$\rho_i(G_{e_i})$, $\rho_{i+1}(G_{e_{i+1}})$ commute. It is elementary to verify that for all  
$g\in G_{e_i}, g'\in G_{e_{i+1}}$ the vectors $\psi_i(g(f_i))$ and $\psi_{i+1}(g(f_{i+1}))$ 
are mutually orthogonal (after being translated to $T_{x_{i+1}}\H^p$ along $e_i, e_{i+1}$). 
The group action $\rho_i(G_{e_i})\acts T_{x_{i+1}}\H^p$ permutes the vectors 
$$
\{\psi_i(g(f_i)), g\in G_{e_i}\}
$$
and fixes the orthogonal complement to these vectors; same is true for the action of 
$G_{e_{i+1}}$ and the vectors 
$$
\{\psi_{i+1}(g'(f_{i+1})), g'\in G_{e_{i+1}}\}
$$
Hence the groups $\rho_i(G_{e_i}), \rho_{i+1}(G_{e_{i+1}})$ commute. 
Therefore we have constructed a homomorphism 
$\rho\co  G\to \isom(\H^p)$, $\rho|G_{e_i}= \rho_i|G_{e_i}$. 

This action has the same ``orthogonality'' properties as 
the homomorphism $\rho$ in the proof of Theorem 
\ref{main}, ie, if $i\ne j$ then for all 
$g\in G_{e_i}\setminus \{1\}, g'\in G_{e_j}\setminus \{1\}$, 
the hyperbolic planes $\H^2, \rho(g)\H^2$ and $\rho(g')\H^2$ are mutually 
orthogonal. Thus, the arguments of the second part of the proof of Theorem \ref{main} 
still work and, by applying Theorem \ref{qi}, 
 we conclude that $\rho$ is discrete, faithful, convex-cocompact. 
\qed 

\medskip 
Suppose now that $X$ is a (locally finite) right-angled 2--dimensional hyperbolic building 
whose fundamental chamber $F$ has $n\ge 6$ vertices. Recall that $X$ is uniquely determined 
by the {\em thickness} $t_i$ of the edges $e_i$ of $F$, ie, the number of 2--faces in $X$ 
containing $e_i$. Thus every such building is the universal cover of an $n$--gon $\p$ 
of finite groups corresponding to a cyclic graph-product. Thickness of the 
edge $e_i$ is the order of the edge group $G_{e_i}$ in $\p$. 

According to a recent theorem of F.\ Haglund, \cite{Haglund}, 
all uniform lattices in the building $X$ are commensurable. 
Hence, as an application of Theorem \ref{graph}, we obtain

\begin{cor}
\label{maincor}
Let $H$ be a group acting discretely, cocompactly and isometrically on $X$. 
Then $H$ contains a finite index subgroup which admits a properly discontinuous  
convex-cocompact action on $\H^p$ for some $p=p(X)$. 
\end{cor}

\section{Extension of discrete representations}
\label{extension}

In this section we discuss the following question:

{\em Suppose that $G\acts \H^n$ is a properly discontinuous isometric action. 
Is it true that $G$ is isomorphic to a Kleinian group?}

\medskip 
Note that the kernel $F$ of the action $G\acts \H^n$ is necessarily finite, therefore 
we have a short exact sequence
$$
1\to F\to G\to \bar{G}\to 1,
$$
where $\bar{G}$ is Kleinian. What we are interested in is whether the group $G$ is itself isomorphic to 
a Kleinian group. 
Of course, a necessary condition for this is that $G$ is residually finite. Finding a 
non-residually finite extension $G$ of a Kleinian group $\bar{G}$ is a very difficult task, 
and presently such extensions are not known. Nevertheless we have:

\begin{thm}
\label{ext}
Suppose that $G$ is a residually finite group which fits into a short exact
sequence
$$
1\to F\to G\to \bar{G}\to 1,
$$
where $\bar{G}$ admits a discrete and faithful representation $\bar\rho$ 
into $\isom(\H^n)$. Then $G$ also admits a discrete and faithful 
representation $\rho$ 
into $\isom(\H^m)$ for some $m$. Moreover, if $\bar\rho$ is convex-cocompact 
(resp.\ geometrically finite) then $\rho$ can be taken convex-cocompact 
(resp.\ geometrically finite). 
\end{thm}
\proof The proof of this theorem is modeled on the proof of the well-known fact that a finite extension 
of a residually finite linear group is again linear, but we present it here for the sake of completeness. 

We first lift $\bar\rho$ to a homomorphism $\bar\rho\co  G\to \isom(\H^n)$, so 
that $Ker(\bar\rho)=F$. 
Since $G$ is residually finite, there exists a homomorphism
$$
\phi\co  G\to Q
$$
where $Q$ is a finite group, so that $\phi|F$ is injective. Embed $Q$ in $SO(k)$ for some $k$. 
The product group $\isom(\H^n)\times SO(k)$ embeds in $\isom(\H^{n+k})$ as the stabilizer of 
$\H^n$ embedded in $\H^{n+k}$ as a totally-geodesic subspace. Therefore, for $m=n+k$ we get a 
homomorphism
$$
\rho\co  G\to \isom(\H^n)\times SO(k)\subset \isom(\H^{m})
$$ 
given by
$$
\rho(g)=(\bar\rho(g), \phi(g)). 
$$
It is clear that $\rho(g)|\H^n= \bar\rho(g)$ and therefore $\rho$ is faithful and  
 $\rho(G)\subset \isom(\H^m)$ is discrete. Moreover, 
$\La(\rho(G))=\La(\bar\rho(\bar{G}))$. Recall that geometrically finite and convex-cocompact 
actions can be detected by considering the dynamics of a discrete group on its limit set 
(see \cite{Bowditch(1993b)}). Therefore, if $\bar\rho$ is convex-cocompact 
(resp.\ geometrically finite) then $\rho$ is also convex-cocompact 
(resp.\ geometrically finite). \qed 

\medskip 
Combining Theorem \ref{ext} with Theorem \ref{main1} we get Theorem \ref{main}. 

\section{Example of a nonlinear Gromov-hyperbolic group}
\label{appendix}

\begin{thm}
\label{non}
There exists an infinite hyperbolic group $G$ such that each representation of $G$ to $GL(m,\k)$ 
factors through a finite group.  In particular, $G$ is nonlinear. Here $\k$ is an abritrary field. 
\end{thm}
\proof Let $\Ga$ be a uniform lattice in a quaternionic hyperbolic space $\H \H^n, n\ge 2$. Since 
$\H \H^n$ is negatively curved, the group $\Ga$ is hyperbolic. Clearly, the group $\Ga$ is a nonelementary 
hyperbolic group; hence $\Ga$ admits an infinite proper quotient $\Ga\to G$ where $G$ is a hyperbolic group 
(see \cite{Gromov(1987)} or \cite{Ivanov-Olshanskii}). We first consider the case when $\k$ has zero characteristic. 
Then without loss of generality we can assume that we are given a linear representation 
$\rho\co  G \to GL(N, \R)$. We will show that $\rho(G)$ is finite 
by using the standard ``adelic'' trick. The reader can find similar applications of this argument 
in Margulis' proof of arithmeticity of higher rank lattices (see \cite{Margulis,Zimmer(1984)}), and 
in Tits' proof of the {\em Tits alternative}, \cite{Tits(1972)}. 

The representation $\rho$ lifts 
to a linear representation $\tilde\rho\co  \Ga\to  GL(N, \R)$. Let $L$ denote 
the Zariski closure of $\tilde\rho(\Ga)$ in $GL(N, \R)$. Let $S$ denote the solvable radical of $L$.     
We first consider the case when $L':=L/S$ is a reductive 
group with nontrivial noncompact factor $H$. Then the projection $\Ga \to G\to L\to H$ 
has Zariski dense image.  Hence,  
according to Corlette's Archimedean superrigidity theorem \cite{Corlette}, 
the representation $\Ga\to H$ extends to $\isom(\H\H^n)$. This however  
contradicts the assumption that the projection $\Ga\to G$ is not 1--1. 
Therefore  the group $L'$ is a compact algebraic group. 
 
Suppose that the projection $\rho(G)\subset L'$ is infinite.  
As a compact Lie group, $L'$   is isomorphic to a subgroup 
of $O(M)$. Since $G$ satisfies property (T), 
$$
H^1(G, o(M)_{Ad(\rho)})=0, 
$$
where $o(M)$ is the Lie algebra of $O(M)$. Vanishing of the above cohomology group
 implies that the space $Hom(G, O(M))/O(M)$ is finite. 
 Hence, analogously to the 
proof of Theorem 7.67 in \cite{Raghunathan}, $\rho$ is conjugate to a representation $\rho'\co  G\to O(M)$  
for which 
$$
\rho'(G)\subset K(F)\subset O(M, F)\subset GL(M, F),
$$
where $F$ is a number field and $K(\R)$ is the Zariski closure of $\rho'(G)$. 

One would like to replace the representation $\rho$ with another representation 
$\phi$ of the group $G$, whose image is Zariski dense in a certain noncompact algebraic group 
and so that $Ker(\phi)=Ker(\rho)$. The most obvious thing to try is to find an element $\si$ 
of the Galois group $Gal(\C/\Q)$, so that the image of $\phi=\si(\rho)$ is not relatively compact. 
This does not necessarily work. Note however, that the restriction of the norm on $\C$ to $\si(F)$ 
gives rise to an {\em Archimedean valuation} on $F$. The idea of the {\em adelic trick} is 
to use {\em non-Archimedean valuations} $v$ of $F$ together with Archimedean ones. This 
is done by introducing the {\em ring of adeles} of $F$, which is a certain subset of the product  
$$
\prod_{v\in Val(F)} F_v,
$$
where $F_v$ is the completion of $F$ with respect to the valuation $v$. 

Let $\A(F)$ denote the ring of adeles of $F$; then the diagonal embedding 
$F\embed \A(F)$ has {\em discrete image} (see for example \cite{Lang}). Hence the diagonal embedding 
$$
\rho'(G)\embed GL(M, \A(F))
$$ 
also has discrete image. If the projection of $\rho'(G)$ to each factor $GL(M, F_v)$ 
were relatively compact, the image of $\rho'(G)$ in $ GL(M, \A(F))$ would be compact as well. 
However a discrete subset of a compact is finite, which contradicts the assumption that $\rho(G)$ 
is infinite. 

Thus there exists a valuation $v$ of $F$ so 
that the image of the projection 
$$
\rho'(G)\to K(F_v)\subset GL(M, F_v)
$$  
is not relatively compact. In case when $v$ is an Archimedean valuation, we can again apply 
Corlette's Archimedean superrigidity \cite{Corlette} to get a contradiction. Hence 
such $v$ has to be nonarchimedean. Therefore the representation $\Ga\to \rho'(G)\to GL(M, F_v)$ 
corresponds to an isometric action of $\Ga$ on a locally finite Euclidean building $X$. However, by the 
non-Archimedean superrigidity theorem of Gromov and Schoen \cite{Gromov-Schoen}, $\Ga$ fixes a point 
in $X$. Therefore the image of $\rho'(G)$ in $GL(M, F_v)$ is relatively compact, which is a contradiction. 
Hence $\rho'(G)$ is finite. It follows that the group $L$ is commensurable to its solvable radical $S$; 
hence $\rho(G)$ is a virtually solvable group. By applying property (T) again, we conclude that 
$\rho(G)$ is finite. 

\medskip
We now consider the case when $\k$ has positive characteristic; since the argument  is similar to the zero characteristic case, 
we give only a sketch. Under the above assumptions, $\rho(G)\subset GL(m, F)$, where $F$ is a finitely generated field 
(of positive characteristic). The field $F$ is an extension
$$
F_q\subset E\subset F
$$ 
where $F_q$ is a finite field, $F_q\subset E$ is an purely transcendental extension 
and $E\subset F$ is an algebraic extension (see \cite[Chapter VI.1]{Hungerford}). 
Since $F$ is finitely generated, $F/E$ is finite-dimensional and therefore, by 
passing to a bigger matrix group, we reduce the problem to the case when $F=E$ is a 
purely transcendental extension, which necessarily has finite transcendence degree. Therefore 
we reduced to the case of    $\rho\co  G\to GL(N, F)$, where $F=F_q(t_1,\ldots,t_m)$ is the field of rational functions 
with coefficients in $F_q$.   

Then we associate with each variable $t_j^{\pm 1}$ a discrete valuation 
$v_{\pm j}$ and an action $G\acts X_{\pm j}$ on the corresponding Euclidean building. The 
non-Archimedean superrigidity theorem of Gromov and Schoen \cite{Gromov-Schoen} shows that 
for each $\pm j$ the action $G\acts X_{\pm j}$ has a fixed point. Therefore the matrix 
coefficients of $\rho(G)$ have bounded degree with respect to all the variables 
$t_1^{\pm 1},\ldots,t_m^{\pm 1}$. Hence, since $F_q$ is finite, the matrix 
coefficients of $\rho(G)$ belongs to a finite subset of $F$ and thus $\rho(G)$ 
is finite. \qed

\begin{rem}
 After this paper was written, I was informed by Alain Valette that he also knew 
how to prove Theorem \ref{non}. I am sure that other people were also aware of this proof 
since all the arguments here are quite standard. 
\end{rem}

\end{document}

%% file: gtoutput.tex
\def\ifplaintex{\expandafter\ifx\csname documentclass\endcsname\relax}

\ifplaintex 
\else
\fi

\expandafter\ifx\csname beginpicture\endcsname\relax
\expandafter\ifx\csname documentclass\endcsname\relax
\input pictex \else
\input prepictex \input pictex \input postpictex \fi\fi

\def\gt{{\mathsurround=0pt\it $\cal G\mskip-2mu$eometry \&\ 
$\cal T\!\!$opology}}        %  journal title in recommended style

\def\gtp{{\mathsurround=0pt\it $\cal G\mskip-2mu$eometry \&\ 
$\cal T\!\!$opology $\cal P\!$ublications}}  % GT publications

\def\lognumber#1{\def\thelognumber{#1}}
\def\volumenumber#1{\def\thevolumenumber{#1}}
\def\papernumber#1{\def\thepapernumber{#1}}
\def\volumeyear#1{\def\thevolumeyear{#1}}

\def\pagenumbers#1#2{\def\startpage{#1}\def\finishpage{#2}}
\def\published#1{\def\publishdate{#1}}
\def\proposed#1{\def\theproposer{#1}}
\def\seconded#1{\def\theseconders{#1}}
\def\received#1{\def\receiveddate{#1}}
\def\revised#1{\def\reviseddate{#1}}
\def\accepted#1{\def\accepteddate{#1}}

\let\\\par\let\thelognumber\relax
\let\thevolumenumber\relax\let\thepapernumber\relax
\let\thevolumeyear\relax\let\thesamplenumber\relax\let\startpage\relax
\let\finishpage\relax\let\publishdate\relax\let\receiveddate\relax
\let\reviseddate\relax\let\accepteddate\relax\let\theasciititle\relax
\let\theasciiauthors\relax
\let\theasciiabstract\relax
\let\theasciiemail\relax\let\theshortauthors\relax\let\theshorttitle\relax

\long\def\maketitlep{   % start of definition of \maketitlep

\count0=\startpage

\gt\hfill      %   Journal title (top left) 
\beginpicture
\setcoordinatesystem units <0.33truein, 0.33truein> point at 2.2 0.9
\setplotsymbol ({$\cal G$})
\plotsymbolspacing=9truept
\circulararc 315 degrees from 0 1 center at 0 0
\setplotsymbol ({$\cal T$})
\circulararc 315 degrees from 1 -1 center at 1 0
\endpicture
\break
{\small\ifx\thesamplenumber\relax % sample?  
Volume \else Sample
\fi\thevolumenumber\ (\thevolumeyear)
\startpage--\finishpage\nl
Published: \publishdate}
\vglue 0.5truein plus 0.4fil minus 0.1truein

{\parskip=0pt\leftskip 0pt plus 1fil\def\\{\par\smallskip}{\ifplaintex\large
\else\Large\fi\bf\thetitle}\par\medskip}   

\vglue 0pt plus 0.1fil 

{\parskip=0pt\leftskip 0pt plus 1fil\def\\{\par}{\sc\theauthors}
\par\medskip}

\vglue 0pt plus 0.1fil 

{\small\parskip=0pt\let\newline\\
{\leftskip 0pt plus 1fil\def\\{\par}{\sl\theaddress}\par}
\expandafter\ifx\theemail\relax    % email address?
\relax\else\vglue 5pt plus 0.02fil minus 2pt\def\\{\stdspace{\rm 
and}\stdspace} 
\cl{Email:\stdspace\tt\theemail}\fi
\ifx\theurl\relax                  % URL given?
\relax\else\vglue 5pt plus 0.02fil minus 2pt\def\\{\stdspace{\rm 
and}\stdspace}
\cl{URL:\stdspace\tt\theurl}\fi\par}

\vglue 7pt plus 0.3fil minus 3pt

{\bf Abstract}
\vglue 5pt plus 0.1fil minus 2pt

\theabstract

\vglue 7pt plus 0.3fil minus 3pt

{\bf AMS Classification numbers}\quad Primary:\quad \theprimaryclass

Secondary:\quad \thesecondaryclass

\vglue 5pt plus 0.3fil minus 2pt

{\bf Keywords:}\quad \thekeywords

\vglue 10pt plus 0.5fil minus 5pt

{\small  Proposed: \theproposer\hfill Received: \receiveddate\nl
Seconded: \theseconders\hfill 
\ifx\reviseddate\relax                         % paper revised?
Accepted: \accepteddate                        % no
\else
Revised: \reviseddate                          % yes
\fi}
\eject
}       %  end of definition of \maketitlep

\let\maketitlepage\maketitlep
\let\maketitle\maketitlepage

\font\phead=cmsl9 scaled 950
\font=cmsl9 scaled 1050
\font\pnum=cmbx10 scaled 913
\font\lnum=cmbx10 
\font\pfoot=cmsl9 scaled 950
\font=cmsl9 scaled 1050
\ifplaintex
\headline{\vbox to 0pt{\vskip -4.5mm\line{\small\phead\ifnum
\count0=\startpage ISSN 1364-0380 (on line)
1465-3060 (printed) \hfill {\pnum\folio}\else\ifodd\count0\def\\{ }% 
\ifx\theshorttitle\relax\thetitle\else\theshorttitle\fi\hfill{\pnum\folio}
\else\def\\{ and }{\pnum\folio}\hfill\ifx\theshortauthors\relax\theauthors
\else\theshortauthors\fi\fi\fi}\vss}}
\footline{\vbox to 0pt{\vglue 0mm\line{\small\pfoot\ifnum\count0=\startpage
\copyright\ \gtp\hfill\else
\gt, Volume \thevolumenumber\ (\thevolumeyear)\hfill\fi}\vss
}}
\else
\makeatletter
\def\@oddhead{{\small\lhead\ifnum\count0=\startpage ISSN 1364-0380 (on line)
1465-3060 (printed) \hfill {\lnum\number\count0}\else\ifodd\count0
\def\\{ }\ifx\theshorttitle\relax \thetitle \else\theshorttitle\fi\hfill
{\lnum\number\count0}\else\def\\{ and }{\lnum\number\count0}
\hfill\ifx\theshortauthors\relax 
\theauthors\else\theshortauthors\fi\fi\fi}}\def\@evenhead{\@oddhead}
\def\@oddfoot{\small\lfoot\ifnum\count0=\startpage\copyright\ \gtp\hfill\else
\gt, Volume \thevolumenumber\ (\thevolumeyear)\hfill\fi}
\def\@evenfoot{\@oddfoot}
\makeatother
\fi

\newwrite\gtoutfile
\long\gdef\makeheadfile{  %%% start of definition of \makeheadfile
{\def\\{, }\def\s{ }
\immediate\openout\gtoutfile head.xxx
\immediate\write\gtoutfile{Proxy-for: \ifx\theasciiauthors\relax
\theauthors\else\theasciiauthors\fi\s<\ifx\theasciiemail\relax\theemail\else\theasciiemail\fi>}
\immediate\write\gtoutfile{\noexpand\\}
\immediate\write\gtoutfile{Authors: \ifx\theasciiauthors\relax
\theauthors\else\theasciiauthors\fi}
{\def\\{ }\immediate\write\gtoutfile{Title: \ifx\theasciititle\relax
\thetitle\else\theasciititle\fi}}
\immediate\write\gtoutfile{Subj-class: GT or SG or MG etc}
\immediate\write\gtoutfile{MSC-class: \theprimaryclass\ifx\thesecondaryclass\relax\else, \thesecondaryclass\fi}
\immediate\write\gtoutfile{Journal-ref: Geom. Topol. \thevolumenumber
(\thevolumeyear) \startpage-\finishpage}
\immediate\write\gtoutfile{Comments: Published by Geometry and Topology at}
\immediate\write\gtoutfile{\s\s http://www.maths.warwick.ac.uk/gt/GTVol\thevolumenumber/paper\thepapernumber.abs.html}
\immediate\write\gtoutfile{\noexpand\\}
\immediate\write\gtoutfile{}
\ifx\theasciiabstract\relax
\immediate\write\gtoutfile{\theabstract}\else
\immediate\write\gtoutfile{\theasciiabstract}\fi
\immediate\write\gtoutfile{}
\immediate\write\gtoutfile{\noexpand\\}
\immediate\write\gtoutfile{}
\immediate\closeout\gtoutfile}}  %%% end of definition of \makeheadfile

\def\maketitlepage{\maketitlep\makeheadfile}
\let\maketitle\maketitlepage